\documentclass[12pt,a4paper]{article}
\usepackage{}
\usepackage{amsfonts}
\usepackage{amssymb}
\usepackage{mathrsfs}
\usepackage{amsmath}
\usepackage{amsthm}
\usepackage{enumerate}
\usepackage{cite}
\usepackage[all]{xy}
\allowdisplaybreaks
\newtheorem{defn}{Definition}[section]
\newtheorem{thm}[defn]{Theorem}
\newtheorem{lem}[defn]{Lemma}
\newtheorem{prop}[defn]{Proposition}
\newtheorem{cor}[defn]{Corollary}
\newtheorem{eg}[defn]{Example}
\newtheorem{re}[defn]{Remark}

\newcommand\relphantom[1]{\mathrel{\phantom{#1}}}

\newcommand{\bdefn}{\begin{defn}}
\newcommand{\edefn}{\end{defn}}
\newcommand{\bthm}{\begin{thm}}
\newcommand{\ethm}{\end{thm}}
\newcommand{\blem}{\begin{lem}}
\newcommand{\elem}{\end{lem}}
\newcommand{\bprop}{\begin{prop}}
\newcommand{\eprop}{\end{prop}}
\newcommand{\bcor}{\begin{cor}}
\newcommand{\ecor}{\end{cor}}
\newcommand{\beg}{\begin{eg}}
\newcommand{\eeg}{\end{eg}}
\newcommand{\bre}{\begin{re}}
\newcommand{\ere}{\end{re}}
\newcommand{\bpf}{\begin{proof}}
\newcommand{\epf}{\end{proof}}

\newcommand{\ch}{{\rm ch}}
\newcommand{\id}{{\rm id}}
\newcommand{\ad}{{\rm ad}}
\newcommand{\End}{{\rm End}}
\newcommand{\K}{\mathbb{K}}

\newcommand{\Z}{\mathbb{Z}}
\newcommand{\g}{\mathfrak{g}}
\newcommand{\A}{\mathscr{A}}
\newcommand{\B}{\mathscr{B}}
\newcommand{\x}{\mathscr{X}}
\newcommand{\y}{\mathscr{Y}}
\newcommand{\z}{\mathscr{Z}}
\newcommand{\F}{\mathscr{F}}
\newcommand{\G}{\mathscr{G}}

\newcommand{\benu}{\begin{enumerate}}
\newcommand{\eenu}{\end{enumerate}}
\newcommand{\bc}{\begin{center}}
\newcommand{\ec}{\end{center}}
\newcommand{\bea}{\begin{eqnarray}}
\newcommand{\eea}{\end{eqnarray}}
\newcommand{\Bea}{\begin{eqnarray*}}
\newcommand{\Eea}{\end{eqnarray*}}
\newcommand{\beq}{\begin{equation}}
\newcommand{\eeq}{\end{equation}}
\newcommand{\Beq}{\begin{equation*}}
\newcommand{\Eeq}{\end{equation*}}
\newcommand{\bspl}{\begin{split}}
\newcommand{\espl}{\end{split}}

\numberwithin{equation}{section}
\newcommand{\supercite}[1]{\textsuperscript{\cite{#1}}}

\begin{document}
\title{{\bf  On the cohomology and extensions of $n$-ary multiplicative Hom-Nambu-Lie superalgebras}}
\author{ Baoling Guan$^{1,2},$  Liangyun Chen$^{1},$   Yao Ma$^{1}$
 \date{{\small {$^1$ School of Mathematics and Statistics, Northeast Normal
 University,\\
Changchun 130024, China}\\{\small {$^2$ College of Sciences, Qiqihar
University, Qiqihar 161006, China}}}}}

\maketitle
\date{}

\begin{abstract}

 In this paper, we discuss the representations of $n$-ary multiplicative Hom-Nambu-Lie superalgebras as a generalization of the notion of representations for $n$-ary multiplicative Hom-Nambu-Lie algebras. We also give the cohomology of an $n$-ary multiplicative Hom-Nambu-Lie superalgebra and obtain a relation between extensions of an $n$-ary multiplicative Hom-Nambu-Lie superalgebra $\mathfrak{b}$ by an abelian one $\mathfrak{a}$ and $Z^1(\mathfrak{b}, \mathfrak{a})_{\bar{0}}$. We also introduce the notion of $T^*$-extensions of $n$-ary multiplicative Hom-Nambu-Lie superalgebras and prove that every finite-dimensional nilpotent metric $n$-ary multiplicative Hom-Nambu-Lie superalgebra $(\g,[\cdot,\cdots,\cdot]_{\g},\alpha,\langle ,\rangle_{\g})$ over  an algebraically closed field of characteristic not 2 in the case $\alpha$ is a surjection is isometric to a suitable $T^*$-extension.
\bigskip

\end{abstract}
\renewcommand{\thefootnote}{\fnsymbol{footnote}}
\footnote[0]{ Corresponding author(L. Chen): chenly640@nenu.edu.cn.}
\footnote[0]{Supported by  NNSF of China (No.11171055),  NSF of  Jilin province (No.201115006), Scientific
Research Foundation for Returned Scholars
    Ministry of Education of China  and the Fundamental Research Funds for the Central Universities (No.12SSXT139). }

\section{Introduction}

In 1996, the concept of $n$-Lie
superalgebras was firstly introduced by Y. Daletskii and V. Kushnirevich in {\rm
  \cite{dykv}}. Moreover,  N. Cantarini and V. G. Kac gave a more general concept of $n$-Lie
superalgebras again  in 2010 in {\cite{CK}}. $n$-Lie superalgebras are more general structures
including $n$-Lie algebras ($n$-ary Nambu-Lie algebras), $n$-ary Nambu-Lie superalgebras and Lie superalgebras.

 The general Hom-algebra structures arose first in connection with quasi-deformation and discretizations of Lie algebras
of vector fields. These quasi-deformations lead to quasi-Lie algebras, a generalized Lie algebra structure in which the skewsymmetry
and Jacobi conditions are twisted. Hom-Lie algebras, Hom-associative algebras, Hom-Lie superalgebras, Hom-bialgebras, $n$-ary Hom-Nambu-Lie algebras and quasi-Hom-Lie algebras are discussed in \cite{aaf,ld,bs,masss,FSA,hls,af,mas,ydd,yddd,sy}. Generalizations of $n$-ary algebras of Lie type and associative type by twisting the identities using linear maps have been introduced in \cite{ah}.

Cohomologies are powerful tools in mathematics, which can be applied to algebras and topologies as well as the theory of smooth manifolds or of holomorphic functions. The cohomology of Lie algebras was defined by C. Chevalley and S. Eilenberg in order to give an algebraic construction of the cohomology of the underlying topological spaces of compact Lie groups\supercite{CE}. The cohomology of Lie superalgebras was introduced by M. Scheunert and R. B. Zhang\supercite{SZ} and was used in mathematics and theoretical physics: the theory of cobordisms, invariant differential operators, central extensions and deformations, etc. The theory of cohomology for $n$-ary Hom-Nambu-Lie algebras and $n$-Lie superalgebras can be found in \cite{AI,YL}. This paper generalizes it to $n$-ary multiplicative Hom-Nambu-Lie superalgebras.

The extension is an important way to find a larger algebra and there are many extensions such as double extensions and Kac-Moody extensions, etc. In 1997, Bordemann introduced the notion of $T^*$-extensions of Lie algebras\supercite{B} and proved that every nilpotent finite-dimensional algebra over an algebraically closed field carrying a nondegenerate invariant symmetric bilinear form is a suitable $T^*$-extension. The method of $T^*$-extension was used in \cite{BBM} and was generalized to many other algebras recently\supercite{LZ2,lyc,YL}. This paper researches general extensions and $T^*$-extensions of $n$-ary multiplicative Hom-Nambu-Lie superalgebras.
In addition, the paper also discusses representations of $n$-ary multiplicative Hom-Nambu-Lie superalgebras as a generalization of the notions of representations for $n$-ary multiplicative Hom-Nambu-Lie algebras.

This paper is organized as follows.  In section 2, we give the representation and the cohomology for an $n$-ary multiplicative Hom-Nambu-Lie superalgebra.
 In section 3, we give a one-to-one correspondence between extensions of an $n$-ary multiplicative Hom-Nambu-Lie superalgebras $\mathfrak{b}$
 by an abelian one $\mathfrak{a}$ and $Z^1(\mathfrak{b}, \mathfrak{a})_{\bar{0}}.$
 In section 4, we introduce the notion of $T^*$-extensions of $n$-ary multiplicative Hom-Nambu-Lie superalgebras and prove that every finite-dimensional
 nilpotent metric $n$-ary multiplicative Hom-Nambu-Lie superalgebra $(\g,[\cdot,\cdots,\cdot],\alpha,\langle ,\rangle_{\g})$
 over  an algebraically closed field of characteristic not 2 such that $\alpha(\g)=\g$ is isometric to (a nondegenerate ideal of codimension 1 of)
 a $T^*$-extension of a nilpotent $n$-ary multiplicative Hom-Nambu-Lie superalgebra whose nilpotent length is at most a half of the nilpotent length of $\g$.


\bdefn \supercite{afn}
An $n$-ary Nambu-Lie superalgebra is a pair $(\g, [\cdot,\cdots,\cdot])$ consisting of a $\Z_2$-graded vector space $\g=\g_{\bar{0}}\oplus \g_{\bar{1}}$ and a multilinear mapping $[\cdot,\cdots,\cdot]: \underbrace{\g\times\cdots\times \g}_{n}\rightarrow \g,$ satisfying
\begin{align*}|[x_1,\cdots, x_n]|=&|x_1|+\cdots+|x_n|,\\
[x_1,\cdots, x_i, x_{i+1},\cdots, x_n]=&-(-1)^{|x_i| |x_{i+1}|}[x_1,\cdots, x_{i+1}, x_i,\cdots, x_n],\\
\begin{split}
[x_1,\cdots, x_{n-1}, [y_1,\cdots, y_n]]=&\sum_{i=1}^n(-1)^{(|x_1|+\cdots+|x_{n-1}|)(|y_1|+\cdots+|y_{i-1}|)}\\
\cdot[y_1,\cdots,y_{i-1},&[x_1,\cdots, x_{n-1},y_i], y_{i+1},\cdots,y_n],
\end{split}\end{align*}
where $|x|\in\Z_2$ denotes the degree of a homogeneous element $x\in\g$.
\edefn

\bdefn
An $n$-ary Hom-Nambu-Lie superalgebra is a triple $(\g, [\cdot,\cdots,\cdot],\alpha)$ consisting of a $\Z_2$-graded vector space $\g=\g_{\bar{0}}\oplus \g_{\bar{1}},$ a multilinear mapping $[\cdot,\cdots,\cdot]: \underbrace{\g\times\cdots\times \g}_{n}\rightarrow \g$ and a family $\alpha=(\alpha_{i})_{1\leq i\leq n-1}$ of even linear maps $\alpha_{i}: \g\rightarrow \g,$ satisfying
\begin{align}|[x_1,\cdots, x_n]|=&|x_1|+\cdots+|x_n|,\label{eq:1}\\
[x_1,\cdots, x_i, x_{i+1},\cdots, x_n]=&-(-1)^{|x_i| |x_{i+1}|}[x_1,\cdots, x_{i+1}, x_i,\cdots, x_n],\label{eq:2}\\
\begin{split}
[\alpha_{1}(x_1),\cdots, \alpha_{n-1}(x_{n-1}), [y_1,\cdots, y_n]]=&\sum_{i=1}^n(-1)^{(|x_1|+\cdots+|x_{n-1}|)(|y_1|+\cdots+|y_{i-1}|)}\label{eq:3}\\
\cdot[\alpha_{1}(y_1),\cdots,\alpha_{i-1}(y_{i-1}),&[x_1,\cdots, x_{n-1},y_i],\alpha_{i}(y_{i+1}),\cdots,\alpha_{n-1}(y_n)],
\end{split}\end{align}
where $|x|\in\Z_2$ denotes the degree of a homogeneous element $x\in\g.$

An $n$-ary Hom-Nambu-Lie superalgebra $(\g, [\cdot,\cdots,\cdot],\alpha)$ is multiplicative, if
 $\alpha=\\(\alpha_{i})_{1\leq i\leq n-1}$ with $\alpha_{1}=\cdots=\alpha_{n-1}=\alpha$ and satisfying
$$\alpha[x_{1},\cdots,x_{n}]=[\alpha(x_{1}),\cdots,\alpha(x_{n})], \forall x_{1},x_{2},\cdots,x_{n}\in \g.$$
\edefn
If the $n$-ary Hom-Nambu-Lie superalgebra $(\g, [\cdot,\cdots,\cdot],\alpha)$ is multiplicative, then the equation (\ref{eq:3}) can be read:
\begin{align}
\begin{split}
[\alpha(x_1),\cdots, \alpha(x_{n-1}), [y_1,\cdots, y_n]]=&\sum_{i=1}^n(-1)^{(|x_1|+\cdots+|x_{n-1}|)(|y_1|+\cdots+|y_{i-1}|)}\label{eq:31}\\
\cdot[\alpha(y_1),\cdots,\alpha(y_{i-1}),&[x_1,\cdots, x_{n-1},y_i],\alpha(y_{i+1}),\cdots,\alpha(y_n)].
\end{split}\tag{${1.3}^{'}$}\end{align}

It is clear that $n$-ary Hom-Nambu-Lie algebras and Hom-Lie
superalgebras are particular cases of $n$-ary Hom-Nambu-Lie superalgebras. In the sequel, when the notation ``$|x|$'' appears,
it means that $x$ is a homogeneous element of degree $|x|$.

\bdefn
Let $(\g, [\cdot,\cdots,\cdot],\alpha)$ and $(\g^{'}, [\cdot,\cdots,\cdot]^{'},\alpha^{'})$ be two $n$-ary Hom-Nambu-Lie superalgebras, where $\alpha=(\alpha_{i})_{1\leq i\leq n-1}$ and $\alpha^{'}=(\alpha^{'}_{i})_{1\leq i\leq n-1}.$ A linear map $f: \g\rightarrow \g$ is an $n$-ary Hom-Nambu-Lie superalgebra morphism if it satisfies
$$f[x_{1},\cdots,x_{n}]=[f(x_{1}),\cdots,f(x_{n})]^{'},$$
$$f\circ \alpha_{i}=\alpha_{i}^{'}\circ f, \forall i=1,\cdots,n-1.$$
\edefn

\beg
Let $(\g, [\cdot,\cdots,\cdot])$ be an $n$-ary Nambu-Lie superalgebra and let $\rho: \g\rightarrow \g$ be  an $n$-ary Nambu-Lie superalgebra endomorphism.
Then $(\g, \rho\circ[\cdot,\cdots,\cdot],\rho)$ is an $n$-ary multiplicative Hom-Nambu-Lie superalgebra.
\bpf
Put $[\cdot,\cdots,\cdot]_{\rho}:=\rho\circ[\cdot,\cdots,\cdot].$ Then
\begin{align*}
&\rho[x_{1},\cdots,x_{n}]_{\rho}=\rho(\rho[x_{1},\cdots,x_{n}])\\
=&\rho[\rho(x_{1}),\cdots, \rho(x_{n})]\\
=&[\rho(x_{1}),\cdots, \rho(x_{n})]_{\rho}.
\end{align*}
Moreover, we have
\begin{align*}
&[\rho(x_{1}),\cdots,\rho(x_{n-1}),[y_{1},\cdots,y_{n}]_{\rho}]_{\rho}\\
=&\rho[\rho(x_{1}),\cdots,\rho(x_{n-1}),{\rho}[y_{1},\cdots,y_{n}]]\\
=&\rho[\rho(x_{1}),\cdots,\rho(x_{n-1}),[\rho(y_{1}),\cdots,\rho(y_{n})]]\\
=&\rho^{2}[x_{1},\cdots,x_{n-1},[y_{1},\cdots,y_{n}]]\\
=&\rho^{2}\Big(\sum_{i=1}^{n}(-1)^{(|x_1|+\cdots+|x_{n-1}|)(|y_1|+\cdots+|y_{i-1}|)}[y_{1},\cdots,[x_{1},\cdots,x_{n-1},y_{i}],\cdots,y_{n}]\Big)\\
=&\sum_{i=1}^{n}(-1)^{(|x_1|+\cdots+|x_{n-1}|)(|y_1|+\cdots+|y_{i-1}|)}\rho[\rho(y_{1}),\cdots,
\rho[x_{1},\cdots,x_{n-1},y_{i}],\cdots,\rho(y_{n})]\\
=&\sum_{i=1}^{n}(-1)^{(|x_1|+\cdots+|x_{n-1}|)(|y_1|+\cdots+|y_{i-1}|)}[\rho(y_{1}),\cdots,
[x_{1},\cdots,x_{n-1},y_{i}]_{\rho},\cdots,\rho(y_{n})]_{\rho}.
\end{align*}
Therefore, $(\g, \rho\circ[\cdot,\cdots,\cdot],\rho)$ is an $n$-ary multiplicative Hom-Nambu-Lie superalgebra.
\epf
\eeg
\bdefn Let $(\g,[\cdot,\cdots,\cdot]_{\g},\alpha)$ be an $n$-ary Hom-Nambu-Lie superalgebra.  A graded subspace $H\subseteq \g$ is a Hom-subalgebra of $(\g,[\cdot,\cdots,\cdot]_{\g},\alpha)$ if $\alpha(H)\subseteq H$
and $H$ is closed under the bracket operation $[\cdot,\cdots,\cdot]_{\g},$ i.e., $[u_{1},\cdots,u_{n}]_{\g}\in H, \forall u_{1},\cdots,u_{n}\in H.$\\
A graded subspace $H\subseteq \g$ is a Hom-ideal of $(\g,[\cdot,\cdots,\cdot]_{\g},\alpha)$ if $\alpha(H)\subseteq H$
and $[u_{1},u_{2},\cdots,\\u_{n}]_{\g}\in H, \forall u_{1}\in H, u_{2},\cdots,u_{n}\in \g.$
\edefn

\section{Cohomology for $n$-ary multiplicative Hom-Nambu-Lie superalgebras}

\bdefn
Let $(\g, [\cdot,\cdots,\cdot],\alpha)$ be an $n$-ary multiplicative Hom-Nambu-Lie superalgebra. $\x=x_1\wedge\cdots\wedge x_{n-1}\in \g^{\wedge^{n-1}}$ is called a fundamental object of $\g$ and $\forall z\in\g, \x\cdot z:=[x_1,\cdots, x_{n-1},z]$. It is clear that $|\x|=|x_1|+\cdots+|x_{n-1}|$.
\edefn

Let $\x=x_1\wedge\cdots\wedge x_{n-1}$ and $\y=y_1\wedge\cdots\wedge y_{n-1}$ be two fundamental objects of $\g$. A bilinear map $[\cdot,\cdot]_{\alpha}: \g^{\wedge^{n-1}}\times \g^{\wedge^{n-1}}\longrightarrow \g^{\wedge^{n-1}}$ defined by
$$[\x,\y]_{\alpha}=\sum_{i=1}^{n-1}(-1)^{|\x|(|y_1|+\cdots+|y_{i-1}|)}
\alpha(y_{1})\wedge\cdots\wedge\alpha(y_{i-1})\wedge
\x \cdot y_{i}\wedge\alpha(y_{i+1})\wedge\cdots\wedge
\alpha(y_{n-1}).$$
 A linear map $\alpha: \g^{\wedge^{n-1}}\longrightarrow \g^{\wedge^{n-1}}$ defined by $\alpha(\x)=\alpha(x_{1})\wedge\cdots\wedge\alpha(x_{n-1}).$
Then $\alpha[\x,\y]_{\alpha}=[\alpha(\x),\alpha(\y)]_{\alpha}$.

\bprop\label{prop:1}Let $(\g, [\cdot,\cdots,\cdot],\alpha)$ be an $n$-ary multiplicative Hom-Nambu-Lie superalgebra.
Suppose that $\x=x_1\wedge\cdots\wedge x_{n-1}$, $\y=y_1\wedge\cdots\wedge y_{n-1}$ and $\z=z_1\wedge\cdots\wedge z_{n-1}$ are fundamental objects of $\g$ and $z$ is an arbitrary element in $\g$. Then \begin{align}
\alpha(\x)\cdot(\y\cdot z)&=(-1)^{|\x||\y|}\alpha(\y)\cdot(\x\cdot z)+[\x,\y]_{\alpha}\cdot \alpha(z),\label{eq:4}\\
[\alpha(\x),[\y, \z]_{\alpha}]_{\alpha}&=(-1)^{|\x||\y|}[\alpha(\y),[\x, \z]_{\alpha}]_{\alpha}+[[\x,\y]_{\alpha}, \alpha(\z)]_{\alpha},\label{eq:5}\\
[\x,\y]_{\alpha}\cdot \alpha(z)&=-(-1)^{|\x||\y|}[\y,\x]_{\alpha}\cdot \alpha(z).\label{eq:6}
\end{align}
\eprop
\bpf It's easy to see that (\ref{eq:4}) is equivalent to (\ref{eq:31}). Note that
\begin{align}
\begin{split}&[\alpha(\x),[\y,\z]_{\alpha}]_{\alpha}=[\alpha(\x),\sum_{i=1}^{n-1}(-1)^{|\y|(|z_1|+\cdots+|z_{i-1}|)}\alpha(z_1)\wedge\cdots\wedge\y\cdot z_i\wedge\\
&\quad \cdots\wedge \alpha(z_{n-1})]_{\alpha}\end{split}\notag\\
=\begin{split}&\sum_{i=1\atop i< j}^{n-1}(-1)^{|\y|(|z_1|+\cdots+|z_{i-1}|)}(-1)^{|\x|(|z_1|+\cdots+|z_{j-1}|+|\y|)}\label{eq:5a}\\
&\quad\cdot\alpha^{2}(z_1)\wedge\!\cdots\!\wedge\alpha(\y\cdot z_i)\wedge\!\cdots\!\wedge \alpha(\x)\cdot \alpha(z_j)\wedge\!\cdots\!\wedge \alpha^{2}(z_{n-1})
\end{split}\\
\begin{split}&+\sum_{i=1\atop j< i}^{n-1}(-1)^{|\y|(|z_1|+\cdots+|z_{i-1}|)}(-1)^{|\x|(|z_1|+\cdots+|z_{j-1}|)}\label{eq:5b}\\
&\quad\cdot\alpha^{2}(z_1)\wedge\!\cdots\!\wedge\alpha(\x)\cdot \alpha(z_j)\wedge\!\cdots\!\wedge\alpha(\y\cdot z_i)\wedge\!\cdots\!\wedge \alpha^{2}(z_{n-1})
\end{split}\\
&+\sum_{i=1}^{n-1}(-1)^{(|\x|+|\y|)(|z_1|+\cdots+|z_{i-1}|)}\alpha^{2}(z_1)\wedge\cdots\wedge\alpha(\x)\cdot(\y\cdot z_i)\wedge\cdots\wedge \alpha^{2}(z_{n-1})\label{eq:5c}.
\end{align}

Similarly,
\begin{align}
&[\alpha(\y),[\x,\z]_{\alpha}]_{\alpha}\notag\\
=\begin{split}&\sum_{i=1\atop i< j}^{n-1}(-1)^{|\x|(|z_1|+\cdots+|z_{i-1}|)}(-1)^{|\y|(|z_1|+\cdots+|z_{j-1}|+|\x|)}\label{eq:5d}\\
&\quad\cdot\alpha^{2}(z_1)\wedge\!\cdots\!\wedge\alpha(\x\cdot z_i)\wedge\!\cdots\!\wedge \alpha(\y)\cdot \alpha(z_j)\wedge\!\cdots\!\wedge \alpha^{2}(z_{n-1})
\end{split}\\
\begin{split}&+\sum_{i=1\atop j< i}^{n-1}(-1)^{|\x|(|z_1|+\cdots+|z_{i-1}|)}(-1)^{|\y|(|z_1|+\cdots+|z_{j-1}|)}\label{eq:5e}\\
&\quad\cdot\alpha^{2}(z_1)\wedge\!\cdots\!\wedge\alpha(\y)\cdot \alpha(z_j)\wedge\!\cdots\!\wedge\alpha(\x\cdot z_i)\wedge\!\cdots\!\wedge \alpha^{2}(z_{n-1})
\end{split}\\
&+\sum_{i=1}^{n-1}(-1)^{(|\x|+|\y|)(|z_1|+\cdots+|z_{i-1}|)}\alpha^{2}(z_1)\wedge\cdots\wedge\alpha(\y)\cdot(\x\cdot z_i)\wedge\cdots\wedge \alpha^{2}(z_{n-1})\label{eq:5f}
\end{align}
and
\beq
[[\x,\y]_{\alpha},\alpha(\z)]_{\alpha}=\sum_{i=1}^{n-1}(-1)^{(|\x|+|\y|)(|z_1|+\cdots+|z_{i-1}|)}\alpha^{2}(z_1)\wedge\cdots\wedge[\x,\y]_{\alpha}\cdot \alpha(z_i)\wedge\cdots\wedge \alpha^{2}(z_{n-1})\label{eq:5g}.
\eeq

It can be checked that (\ref{eq:5a})+(\ref{eq:5b})=$(-1)^{|\x||\y|}$(\ref{eq:5d})+$(-1)^{|\x||\y|}$(\ref{eq:5e}). By (\ref{eq:4}), we conclude (\ref{eq:5c})=(\ref{eq:5g})+$(-1)^{|\x||\y|}$(\ref{eq:5f}). Thus (\ref{eq:5}) holds.

 Using (\ref{eq:4}), by exchanging $\x$ and $\y$, we have
\beq\alpha(\y)\cdot(\x\cdot z)=[\y,\x]_{\alpha}\cdot \alpha(z)+(-1)^{|\x||\y|}\alpha(\x)\cdot(\y\cdot
z).\label{eq:7}\eeq Comparing (\ref{eq:4}) with (\ref{eq:7}), we
obtain (\ref{eq:6}). \epf

\bdefn
Let $(\g, [\cdot,\cdots,\cdot],\alpha)$ be an $n$-ary multiplicative Hom-Nambu-Lie superalgebra and $V=V_{\bar{0}}\oplus V_{\bar{1}}$ be a $\Z_2$-graded vector space over a field $\K$. A graded representation $\rho$ of $\g$ on $V$ is a linear map $\rho: \g^{\wedge^{n-1}}\rightarrow \End(V), \x\mapsto\rho(\x)=\rho(x_1,\cdots,x_{n-1})$ such that
\begin{align}
&\rho(\x)\cdot V_\beta\subseteq V_{\beta+|\x|}, \forall \beta\in \Z_2,\label{eq:rho1}\\
\rho(\alpha(\x))&\rho(\y)=(-1)^{|\x||\y|}\rho(\alpha(\y))\rho(\x)+
\rho[\x,\y]_{\alpha}\circ \nu,\label{eq:rho2}\\
\begin{split}\label{eq:rho3}
\rho(\alpha(x_1),\cdots,\alpha(x_{n-2}),&[y_1,\cdots,y_n])\circ \nu=
\sum_{i=1}^n(-1)^{n-i}(-1)^{(|x_1|+\cdots+|x_{n-2}|)(|y_1|+\cdots+\widehat{|y_i|}+\cdots+|y_n|)}\\
&\cdot(-1)^{|y_i|(|y_{i+1}|+\cdots+|y_n|)}\rho(\alpha(y_1),\!\cdots\!,\widehat{\alpha(y_i)},\!\cdots\!,\alpha(y_n))\rho(x_1,\!\cdots\!,x_{n-2},y_i),\end{split}
\end{align}
for  $\x,\y\in\g^{\wedge^{n-1}}$ and
$x_1,\cdots,x_{n-2},y_1,\cdots,y_n\in\g,$ $\nu\in \mathrm{End}V,$ where the sign
\textasciicircum ~ indicates that the element below must be
omitted. The $\Z_2$-graded representation space $(V,\nu)$ is said to be a
graded $\g$-module.
\edefn

We use a supersymmetric notation $[x_1,\cdots,x_{n-1}, v]$(like
(\ref{eq:2})) to denote $\rho(\x)\cdot v$ and set
$[x_1,\cdots,x_{n-2}, v_1,v_2]=0$ and $(\alpha+\nu)(x+v)=\alpha(x)+\nu(v)$ for all $x\in\g$ and $v\in V,$ then $(\g\oplus V, [\cdot,\cdots,\cdot], \alpha+\nu)$ becomes an
$n$-ary multiplicative Hom-Nambu-Lie superalgebra such that $V$ is a $\Z_2$-graded
abelian ideal of $\g$, that is,
$$[V,\underbrace{\g,\cdots,\g}_{n-1}]\subseteq V \quad \text{and}\quad [V,V,\underbrace{\g,\cdots,\g}_{n-2}]=0.$$
In the sequel, we will usually abbreviate $\rho(\x)\cdot v$ with
$\x\cdot v$.
\beg
Let $(\g, [\cdot,\cdots,\cdot], \alpha)$ be an $n$-ary multiplicative Hom-Nambu-Lie superalgebra. Then $\mathrm{ad}: \g^{\wedge^{n-1}}\rightarrow \End(\g), \x\mapsto \mathrm{ad}\x$ defined by $$ \mathrm{ad}\x(z)=\x\cdot z$$ is a graded representation of $\g,$ it is also called the adjoint graded representation of $\g.$
\eeg

\bdefn Let $(\g, [\cdot,\cdots,\cdot], \alpha)$ be an $n$-ary multiplicative Hom-Nambu-Lie superalgebra and $(V,\nu)$ be a
graded $\g$-module.
An $m$-cochain is an $(m+1)$-linear map $$f: \underbrace{\g^{\wedge^{n-1}}\otimes\cdots\otimes\g^{\wedge^{n-1}}}_m\wedge\g\longrightarrow V$$
such that $$\nu\circ f(\x_{1},\x_{2},\cdots,\x_{m},z)=f(\alpha(\x_{1}),\alpha(\x_{2}),\cdots,\alpha(\x_{m}),\alpha(z))$$ for all $\x_{1},\x_{2},\cdots,\x_{m}\in \g^{\wedge^{n-1}}$ and $z\in \g.$ We denote the set of $m$-cochain by $C^m(\g, V).$
\edefn

\bdefn\label{def:delta}For $m\geq 1,$ we call $m$-coboundary operator of the $n$-ary multiplicative Hom-Nambu-Lie superalgebra $(\g, [\cdot,\cdots,\cdot],\alpha)$
the even linear map $\delta^{m}: C^m(\g, V)\rightarrow C^{m+1}(\g, V)$ by
\begin{align*}
&(\delta^{m} f)(\x_1,\cdots,\x_m, \x_{m+1}, z)\\
=&\sum_{i<j}(-1)^i(-1)^{|\x_i|(|\x_{i+1}|+\cdots+|\x_{j-1}|)}f(\alpha(\x_1),\cdots,\widehat{\alpha(\x_i)},\cdots,[\x_{i},\x_{j}]_{\alpha},\cdots,\alpha(\x_{m+1}),\alpha(z))\\
&+\sum_{i=1}^{m+1}(-1)^i(-1)^{|\x_i|(|\x_{i+1}|+\cdots+|\x_{m+1}|)}f(\alpha(\x_1),\cdots,\widehat{\alpha(\x_i)},\cdots,\alpha(\x_{m+1}),\x_i\cdot z)\\
& +\sum_{i=1}^{m+1}(-1)^{i+1}(-1)^{|\x_i|(|f|+|\x_{1}|+\cdots+|\x_{i-1}|)}\alpha^{m}(\x_i)\cdot f(\x_1,\cdots,\widehat{\x_i},\cdots,\x_{m+1}, z)\\
&  +(-1)^m(f(\x_1,\cdots,\x_m, ~~)\cdot \x_{m+1})\bullet_{\alpha} \alpha^{m}(z),
\end{align*}
where $\x_i=\x_i^1\wedge\cdots\wedge\x_i^{n-1}\in\g^{\wedge^{n-1}}, i=1,\cdots,m+1, z\in\g$ and the last term is defined by
\Beq\begin{split}
(f(\x_1,\cdots,\x_m,~~ )\cdot \x_{m+1})\bullet_{\alpha} \alpha^{m}(z)=&\sum_{i=1}^{n-1}(-1)^{(|f|+|\x_{1}|+\cdots+|\x_{m}|)(|\x_{m+1}^1|+\cdots+|\x_{m+1}^{i-1}|)}\\
\cdot[\alpha^{m}(\x_{m+1}^1),\cdots,&f(\x_1,\cdots,\x_m,\x_{m+1}^i),\cdots,\alpha^{m}(\x_{m+1}^{n-1}),\alpha^{m}(z)].
\end{split}\Eeq
\edefn

We now check that $\delta^{m+1}\circ \delta^{m}=0$. In fact, for $f\in C^{m}(\g, V),$ we have
\begin{align}
&(\delta^{m+1}\circ \delta^{m}(f))(\x_1,\cdots, \x_{m+2}, z)\notag\\
=&\sum_{i<j}(-1)^i(-1)^{|\x_i|(|\x_{i+1}|+\cdots+|\x_{j-1}|)}\delta^{m} f(\alpha(\x_1),\cdots,\widehat{\alpha(\x_i)},\cdots,[\x_{i},\x_{j}]_{\alpha},\cdots,\notag\\\
&\quad \quad \alpha(\x_{m+2}),\alpha(z))\notag\\
&+\sum_{i=1}^{m+2}(-1)^i(-1)^{|\x_i|(|\x_{i+1}|+\cdots+|\x_{m+2}|)}\delta^{m} f(\alpha(\x_1),\cdots,\widehat{\alpha(\x_i)},\cdots,\alpha(\x_{m+2}),\x_i\cdot z)\notag\\
& +\sum_{i=1}^{m+2}(-1)^{i+1}(-1)^{|\x_i|(|f|+|\x_{1}|+\cdots+|\x_{i-1}|)}\alpha^{m+1}(\x_i)\cdot\delta^{m} f(\x_1,\cdots,\widehat{\x_i},\cdots,\x_{m+2}, z)\notag\\
&  +(-1)^{m+1}(\delta^{m} f(\x_1,\cdots,\x_{m+1}, ~~)\cdot \x_{m+2})\bullet_{\alpha} \alpha^{m+1}(z),\notag\\
\begin{split}
 &=\sum_{s<t<i<j}a_{ijst}f(\alpha^{2}(\x_1),\cdots,\widehat{\alpha^{2}(\x_s)},\cdots,[\alpha(\x_{s}),\alpha(\x_{t})]_{\alpha},\cdots, \widehat{\alpha^{2}(\x_i}),\cdots,\\
 &\relphantom{}\quad \quad   \alpha[\x_{i},\x_{j}]_{\alpha},\cdots, \alpha^{2}(\x_{m+2}),\alpha^{2}(z))\end{split}\tag{a1}\\
  \begin{split}&+\sum_{s<i<t<j}\widetilde{a_{ijst}}f(\alpha^{2}(\x_1),\cdots,\widehat{\alpha^{2}(\x_s)},\cdots, \widehat{\alpha^{2}(\x_i)},\cdots,[\alpha(\x_{s}),\alpha(\x_{t})]_{\alpha}, \cdots,
  \\&\relphantom{}\quad \quad   \alpha[\x_{i},\x_{j}]_{\alpha},\cdots, \alpha^{2}(\x_{m+2}),\alpha^{2}(z))\end{split}\tag{a2}\\
   \begin{split}&+\sum_{s<i<j<t}a_{ijst}f(\alpha^{2}(\x_1),\cdots,\widehat{\alpha^{2}(\x_s)},\cdots, \widehat{\alpha^{2}(\x_i)},\cdots,\alpha[\x_{i},\x_{j}]_{\alpha}, \cdots,\\
 &\relphantom{}\quad \quad [\alpha(\x_{s}),\alpha(\x_{t})]_{\alpha},\cdots, \alpha^{2}(\x_{m+2}),\alpha^{2}(z))\end{split}\tag{a3}\\
    \begin{split}&-\sum_{i<s<t<j}a_{ijst}f(\alpha^{2}(\x_1),\cdots,\widehat{\alpha^{2}(\x_i)},\cdots, \widehat{\alpha^{2}(\x_s)},\cdots,[\alpha(\x_{s}),\alpha(\x_{t})]_{\alpha}, \cdots,\\
    &\relphantom{}\quad \quad  \alpha[\x_{i},\x_{j}]_{\alpha},\cdots,\alpha^{2}(\x_{m+2}),\alpha^{2}(z))\end{split}\tag{a4}\\
      \begin{split}&-\sum_{i<s<j<t}\widetilde{a_{ijst}}f(\alpha^{2}(\x_1),\cdots,\widehat{\alpha^{2}(\x_i)},\cdots, \widehat{\alpha^{2}(\x_s)},\cdots,\alpha[\x_{i},\x_{j}]_{\alpha}, \cdots, \\
    &\relphantom{}\quad \quad[\alpha(\x_{s}),\alpha(\x_{t})]_{\alpha},\cdots,\alpha^{2}(\x_{m+2}),\alpha^{2}(z))\end{split}\tag{a5}\\
      \begin{split}&-\sum_{i<j<s<t}a_{ijst}f(\alpha^{2}(\x_1),\cdots,\widehat{\alpha^{2}(\x_i)},\cdots,\alpha[\x_{i},\x_{j}]_{\alpha},\cdots, \widehat{\alpha^{2}(\x_s)}, \\
    &\relphantom{}\quad \quad \cdots, [\alpha(\x_{s}),\alpha(\x_{t})]_{\alpha},\cdots,\alpha^{2}(\x_{m+2}),\alpha^{2}(z))\end{split}\tag{a6}\\
 \begin{split}&+\sum_{k<i<j}\widetilde{b_{ijk}}f(\alpha^{2}(\x_1),\cdots, \widehat{\alpha^{2}(\x_k)},\cdots, \widehat{\alpha^{2}(\x_i)}, \cdots, [\alpha(\x_k),[\x_{i},\x_{j}]_{\alpha}]_{\alpha}, \\
    &\relphantom{}\quad \quad \cdots,\alpha^{2}(\x_{m+2}),\alpha^{2}(z))\end{split}\tag{b1}\\
   \begin{split}&-\sum_{i<k<j}b_{ijk}f(\alpha^{2}(\x_1),\cdots, \widehat{\alpha^{2}(\x_i)},\cdots, \widehat{\alpha^{2}(\x_k)}, \cdots, [\alpha(\x_k),[\x_{i},\x_{j}]_{\alpha}]_{\alpha},\\
    &\relphantom{}\quad \quad \cdots,\alpha^{2}(\x_{m+2}),\alpha^{2}(z))\end{split}\tag{b2}\\
    \begin{split}&-\sum_{i<j<k}\widetilde{b_{ikj}}f(\alpha^{2}(\x_1),\cdots, \widehat{\alpha^{2}(\x_i)},\cdots, \widehat{\alpha^{2}(\x_j)}, \cdots, [\x_{i},\x_{j}]_{\alpha}\cdot\alpha(\x_k),\\
    &\relphantom{}\quad \quad \cdots, \alpha^{2}(\x_{m+2}),\alpha^{2}(z))\end{split}\tag{b3}\\
  \begin{split}&+\sum_{k<i<j}c_{ijk}f(\alpha^{2}(\x_1),\cdots,\widehat{\alpha^{2}(\x_k)},\cdots,\widehat{\alpha^{2}(\x_i)},\cdots,\alpha[\x_{i},\x_{j}]_{\alpha},\\
    &\relphantom{}\quad \quad\cdots,\alpha^{2}(\x_{m+2}),\alpha(\x_k)\cdot \alpha(z))\end{split}\tag{c1}\\
   \begin{split}&-\sum_{i<k<j}\widetilde{c_{ijk}}f(\alpha^{2}(\x_1),\cdots,\widehat{\alpha^{2}(\x_i)},\cdots,\widehat{\alpha^{2}(\x_k)},\cdots,\alpha[\x_{i},\x_{j}]_{\alpha},\\
   &\relphantom{}\quad \quad\cdots,\alpha^{2}(\x_{m+2}), \alpha(\x_k)\cdot \alpha(z))\end{split}\tag{c2}\\
    \begin{split}&-\sum_{i<j<k}c_{ijk}f(\alpha^{2}(\x_1),\cdots,\widehat{\alpha^{2}(\x_i)},\cdots,\alpha[\x_{i},\x_{j}]_{\alpha}
    ,\cdots,\widehat{\alpha^{2}(\x_k)},\\
    &\relphantom{}\quad \quad\cdots,\alpha^{2}(\x_{m+2}), \alpha(\x_k)\cdot \alpha(z))\end{split}\tag{c3}\\
 \begin{split}&-\sum_{i<j}\widetilde{d_{ij}}f(\alpha^{2}(\x_1),\cdots,\widehat{\alpha^{2}(\x_i)},\cdots,\widehat{\alpha^{2}(\x_j)},\cdots,\alpha^{2}(\x_{m+2}), [\x_{i},\x_{j}]_{\alpha}\cdot \alpha(z))\end{split}\tag{d1}\\
  \begin{split}&+\sum_{k<i<j}e_{ijk}\alpha^{m+1}(\x_k)\cdot f(\alpha(\x_1),\cdots,\widehat{\alpha(\x_k)},\cdots,\widehat{\alpha(\x_i)},\cdots,[\x_{i},\x_{j}]_{\alpha},\\
    &\relphantom{}\quad \quad \cdots,\alpha(\x_{m+2}), \alpha(z))\end{split}\tag{e1}\\
  \begin{split}&-\sum_{i<k<j}\widetilde{e_{ijk}}\alpha^{m+1}(\x_k)\cdot f(\alpha(\x_1),\cdots,\widehat{\alpha(\x_i)},\cdots,\widehat{\alpha(\x_k)},\cdots,[\x_{i},\x_{j}]_{\alpha},\\
    &\relphantom{}\quad \quad \cdots,\alpha(\x_{m+2}), \alpha(z))\end{split}\tag{e2}\\
   \begin{split}&-\sum_{i<j<k}e_{ijk}\alpha^{m+1}(\x_k)\cdot f(\alpha(\x_1),\cdots,\widehat{\alpha(\x_i)},\cdots,[\x_{i},\x_{j}]_{\alpha},\cdots,\widehat{\alpha(\x_k)},\\
    &\relphantom{}\quad \quad\cdots,\alpha(\x_{m+2}), \alpha(z))\end{split}\tag{e3}\\
  \begin{split}&-\sum_{i<j}\widetilde{g_{ij}}\alpha^{m}([\x_{i},\x_{j}]_{\alpha})\cdot f(\alpha(\x_1),\cdots,\widehat{\alpha(\x_i)},\cdots,\widehat{\alpha(\x_j)},\cdots,\alpha(\x_{m+2}), \alpha(z))\end{split}\tag{g1}\\
 \begin{split}&+\sum_{i<j\leq {m+1}}h_{ij}(f(\alpha(\x_1),\cdots,\widehat{\alpha(\x_i)},\cdots,[\x_{i},\x_{j}]_{\alpha},\cdots,\alpha(\x_{m+1}),~~)\cdot\alpha(\x_{m+2}))\\
    &\relphantom{}\quad \quad\bullet_{\alpha} \alpha^{m+1}(z)\end{split}\tag{h1}\\
 \begin{split}
 &+\sum_{k=1}^{m+1}(-1)^{k+m}(-1)^{|\x_k|(|\x_{k+1}|+\cdots+|\x_{m+1}|)}\\
 &\relphantom{+\sum_{k=1}^{m+1}}\cdot(f(\alpha(\x_1),\cdots,\widehat{\alpha(\x_k)},\cdots,\alpha(\x_{m+1}),~~)\cdot[\x_k,\x_{m+2}]_{\alpha})\bullet_{\alpha} \alpha^{m+1}(z)\end{split}\tag{l1}\label{l1}\\
  \begin{split}&+\sum_{s<t<i}c_{sti}f(\alpha^{2}(\x_1),\cdots,\widehat{\alpha^{2}(\x_s)},\cdots,[\alpha(\x_{s}),\alpha(\x_{t})]_{\alpha},
  \cdots,\widehat{\alpha^{2}(\x_i)},\cdots,\\
    &\relphantom{}\quad \quad\alpha^{2}(\x_{m+2}), \alpha(\x_i\cdot z))\end{split}\tag{c4}\\
   \begin{split}&+\sum_{s<i<t}\widetilde{c_{sti}}f(\alpha^{2}(\x_1),\cdots,\widehat{\alpha^{2}(\x_s)},\cdots,\widehat{\alpha^{2}(\x_i)},\cdots,
   [\alpha(\x_{s}),\alpha(\x_{t})]_{\alpha},\cdots,\\
    &\relphantom{}\quad \quad\alpha^{2}(\x_{m+2}), \alpha(\x_i\cdot z))\end{split}\tag{c5}\\
    \begin{split}&-\sum_{i<s<t}c_{sti}f(\alpha^{2}(\x_1),\cdots,\widehat{\alpha^{2}(\x_i)},\cdots,\widehat{\alpha^{2}(\x_s)},\cdots,
    [\alpha(\x_{s}),\alpha(\x_{t})]_{\alpha},\cdots,\\
    &\relphantom{}\quad \quad\alpha^{2}(\x_{m+2}), \alpha(\x_i\cdot z))\end{split}\tag{c6}\\
   \begin{split}&+\sum_{k<i}\widetilde{d_{ik}}f(\alpha^{2}(\x_1),\cdots,\widehat{\alpha^{2}(\x_k)},\cdots,\widehat{\alpha^{2}(\x_i)},\cdots,\alpha^{2}(\x_{m+2}), \alpha(\x_k)\cdot(\x_i\cdot z))\end{split}\tag{d2}\\
 \begin{split}&-\sum_{i<k}d_{ik}f(\alpha^{2}(\x_1),\cdots,\widehat{\alpha^{2}(\x_i)},\cdots,\widehat{\alpha^{2}(\x_k)},\cdots,\alpha^{2}(\x_{m+2}), \alpha(\x_k)\cdot(\x_i\cdot z))\end{split}\tag{d3}\\
  \begin{split}&+\sum_{k<i}p_{ki}\alpha^{m+1}(\x_k)\cdot f(\alpha(\x_1),\cdots,\widehat{\alpha(\x_k)},\cdots,\widehat{\alpha(\x_i)},\cdots,\alpha(\x_{m+2}), \x_i\cdot z)\end{split}\tag{p1}\\
   \begin{split}&-\sum_{i<k}\widetilde{p_{ki}}\alpha^{m+1}(\x_k)\cdot f(\alpha(\x_1),\cdots,\widehat{\alpha(\x_i)},\cdots,\widehat{\alpha(\x_k)},\cdots,\alpha(\x_{m+2}), \x_i\cdot z)\end{split}\tag{p2}\\
\begin{split}&+\sum_{i=1}^{m+1} (-1)^{i+m}(-1)^{|\x_i|(|\x_{i+1}|+\cdots+|\x_{m+2}|)}\\
 &\relphantom{-\sum_{i=1}^{n-1}}\cdot(f(\alpha(\x_1),\cdots,\widehat{\alpha(\x_i)},\cdots,\alpha(\x_{m+1}),~~)\cdot\alpha(\x_{m+2}))
 \bullet_{\alpha}\alpha^{m}(\x_i\cdot z)\end{split}\tag{l2}\label{l2}\\
 &+(f(\alpha(\x_1),\cdots,\alpha(\x_{m}),~~)\cdot \alpha(\x_{m+1}))\bullet_{\alpha}\alpha^{m}(\x_{m+2}\cdot z)\tag{q1}\\
 \begin{split}&+\sum_{s<t<i}e_{sti}\alpha^{m+1}(\x_i)\cdot f(\alpha(\x_1),\cdots,\widehat{\alpha(\x_s)},\cdots, [\x_{s},\x_{t}]_{\alpha},\cdots,\widehat{\alpha(\x_i)},\\
    &\relphantom{}\quad \quad\cdots,\alpha(\x_{m+2}), \alpha(z))\end{split}\tag{e4}\\
   \begin{split}&+\sum_{s<i<t}\widetilde{e_{sti}}\alpha^{m+1}(\x_i)\cdot f(\alpha(\x_1),\cdots,\widehat{\alpha(\x_s)},\cdots,\widehat{\alpha(\x_i)},\cdots, [\x_{s},\x_{t}]_{\alpha}\\
    &\relphantom{}\quad \quad\cdots,\alpha(\x_{m+2}), \alpha(z))\end{split}\tag{e5}\\
   \begin{split}&-\sum_{i<s<t}e_{sti}\alpha^{m+1}(\x_i)\cdot f(\alpha(\x_1),\cdots,\widehat{\alpha(\x_i)},\cdots,\widehat{\alpha(\x_s)},\cdots, [\x_{s},\x_{t}]_{\alpha}\\
    &\relphantom{}\quad \quad\cdots,\alpha(\x_{m+2}), \alpha(z))\end{split}\tag{e6}\\
   &+\sum_{k<i}\widetilde{p_{ik}}\alpha^{m+1}(\x_i)\cdot f(\alpha(\x_1),\cdots,\widehat{\alpha(\x_k)},\cdots,\widehat{\alpha(\x_i)},\cdots,\alpha(\x_{m+2}), \x_k\cdot z)\tag{p3}\\
  &-\sum_{i<k}p_{ik}\alpha^{m+1}(\x_i)\cdot f(\alpha(\x_1),\cdots,\widehat{\alpha(\x_i)},\cdots,\widehat{\alpha(\x_k)},\cdots,\alpha(\x_{m+2}), \x_k\cdot z)\tag{p4}\\
 &+\sum_{k<i}\widetilde{g_{ki}}\alpha^{m+1}(\x_i)\cdot(\alpha^{m}(\x_k)\cdot f(\x_1,\cdots,\widehat{\x_k},\cdots,\widehat{\x_i},\cdots,\x_{m+2},z))\tag{g2}\\
  &-\sum_{i<k}g_{ki}\alpha^{m+1}(\x_i)\cdot(\alpha^{m}(\x_k)\cdot  f(\x_1,\cdots,\widehat{\x_i},\cdots,\widehat{\x_k},\cdots,\x_{m+2},z))\tag{g3}\\
  \begin{split}
 &-\sum_{i=1}^{m+1} (-1)^{i+m}(-1)^{|\x_i|(|f|+|\x_{1}|+\cdots+|\x_{i-1}|)}\\
 &\relphantom{-\sum_{i=1}^{n-1}}\alpha^{m+1}(\x_i)\cdot((f(\x_1,\cdots,\widehat{\x_i},\cdots,\x_{m+1},~~)\cdot \x_{m+2})\bullet_{\alpha} \alpha^{m}(z))
 \end{split}\tag{l3}\label{l3}\\
 &-(-1)^{|\x_{m+2}|(|f|+|\x_1|+\cdots+|\x_{m+1}|)}\alpha^{m+1}(\x_{m+2})\cdot((f(\x_1,\cdots,\x_{m},~~)\cdot \x_{m+1})\bullet_{\alpha} \alpha^{m}(z))\tag{q2}\\
 \begin{split} &-\sum_{s<t\leq {m+1}}h_{st}(f(\alpha(\x_1),\cdots,\widehat{\alpha(\x_s)},\cdots,[\x_{s},\x_{t}]_{\alpha},\cdots,\alpha(\x_{m+1}),~~)\cdot\alpha(\x_{m+2}))\\
    &\relphantom{}\quad \quad\bullet_{\alpha} \alpha^{m+1}(z) \end{split}\tag{h2}\\
 \begin{split}
 &-\!\!\sum_{i=1}^{n-1}\sum_{k=1}^{m+1} (-1)^{m+k}(-1)^{(|f|+|\x_{1}|+\cdots+|\x_{{m+1}}|)(|\x_{m+2}^1|+\cdots+|\x_{m+2}^{i-1}|)} (-1)^{|\x_k|(|\x_{k+1}|+\cdots+|\x_{m+1}|)}\\
 &\relphantom{}\cdot[\alpha^{m+1}(\x_{m+2}^1),\cdots,f(\alpha(\x_1),\cdots,\widehat{\alpha(\x_k)},\cdots,\alpha(\x_{m+1}), \x_k\cdot\x_{m+2}^i),\cdots,\\
 &\relphantom{}\quad \quad\alpha^{m+1}(\x_{m+2}^{n-1}),\alpha^{m+1}(z)]
 \end{split}\tag{l4}\label{l4}\\
 \begin{split}
 &+\!\!\sum_{i=1}^{n-1}\sum_{k=1}^{m+1} (-1)^{m+k}(-1)^{(|f|+|\x_{1}|+\!\cdots+|\x_{{m+1}}|)(|\x_{m+2}^1|+\cdots+|\x_{m+2}^{i-1}|)} (-1)^{|\x_k|(|f|+|\x_1|+\cdots+|\x_{k-1}|)}\\
 &\relphantom{}\cdot[\alpha^{m+1}(\x_{m+2}^1),\cdots,\alpha^{m}(\x_k)\cdot f(\x_1,\cdots,\widehat{\x_k},\cdots,\x_{m+1},\x_{m+2}^i),\cdots,\\
 &\relphantom{}\quad \quad\alpha^{m+1}(\x_{m+2}^{n-1}),\alpha^{m+1}(z)]
 \end{split}\tag{l5}\label{l5}\\
 \begin{split}
 &-\sum_{i=1}^{n-1}(-1)^{(|f|+|\x_{1}|+\cdots+|\x_{{m+1}}|)(|\x_{m+2}^1|+\cdots+|\x_{m+2}^{i-1}|)}\notag\\
 &\relphantom{}\cdot[\alpha^{m+1}(\x_{m+2}^1),\cdots,(f(\x_1,\cdots,\x_{m},~~)\cdot\x_{m+1})\bullet_{\alpha}\alpha^{m}(\x_{m+2}^i),\cdots,\\
 &\relphantom{}\quad \quad
 \alpha^{m+1}(\x_{m+2}^{n-1}),\alpha^{m+1}(z)],
 \end{split}\tag{q3}
\end{align}
where
\begin{align*}
a_{ijst}=&(-1)^{s+i}(-1)^{|\x_i|(|\x_{i+1}|+\cdots+|\x_{j-1}|)}(-1)^{|\x_s|(|\x_{s+1}|+\cdots+|\x_{t-1}|)},
& \widetilde{a_{ijst}}=&(-1)^{|\x_i||\x_s|}a_{ijst};\\
b_{ijk}=&(-1)^{i+k}(-1)^{|\x_i|(|\x_{i+1}|+\cdots+|\x_{j-1}|)}(-1)^{|\x_k|(|\x_{k+1}|+\cdots+|\x_{j-1}|)},
& \widetilde{b_{ijk}}=&(-1)^{|\x_i||\x_k|}b_{ijk};\\
c_{ijk}=&(-1)^{i+k}(-1)^{|\x_i|(|\x_{i+1}|+\cdots+|\x_{j-1}|)}(-1)^{|\x_k|(|\x_{k+1}|+\cdots+|\x_{m+2}|)},
& \widetilde{c_{ijk}}=&(-1)^{|\x_i||\x_k|}c_{ijk};\\
d_{ij}=&(-1)^{i+j}(-1)^{|\x_i|(|\x_{i+1}|+\cdots+|\x_{m+2}|)}(-1)^{|\x_j|(|\x_{j+1}|+\cdots+|\x_{m+2}|)},
& \widetilde{d_{ij}}=&(-1)^{|\x_i||\x_j|}d_{ij};\\
e_{ijk}=&(-1)^{i+k+1}(-1)^{|\x_i|(|\x_{i+1}|+\cdots+|\x_{j-1}|)}(-1)^{|\x_k|(|f|+|\x_{1}|+\cdots+|\x_{k-1}|)},
& \widetilde{e_{ijk}}=&(-1)^{|\x_i||\x_k|}e_{ijk};\\
g_{ij}=&(-1)^{i+j+1}(-1)^{|\x_i|(|f|+|\x_{1}|+\cdots+|\x_{i-1}|)}(-1)^{|\x_j|(|f|+|\x_{1}|+\cdots+|\x_{j-1}|)},
& \widetilde{g_{ij}}=&(-1)^{|\x_i||\x_j|}g_{ij};\\
h_{ij}=&(-1)^{i+m}(-1)^{|\x_i|(|\x_{i+1}|+\cdots+|\x_{j-1}|)},
& \widetilde{h_{ij}}=&(-1)^{|\x_i||\x_j|}h_{ij};\\
p_{ki}=&(-1)^{i+k+1}(-1)^{|\x_i|(|\x_{i+1}|+\cdots+|\x_{m+2}|)}(-1)^{|\x_k|(|f|+|\x_{1}|+\cdots+|\x_{k-1}|)},
&\widetilde{p_{ki}}=&(-1)^{|\x_i||\x_k|}p_{ki}.
\end{align*}

It can be verified that the sum of terms labeled with the same letter vanishes. For example, (\ref{l1})+(\ref{l2})+(\ref{l3})+(\ref{l4})+(\ref{l5})=0,
in fact, \begin{align}
&(\mathrm{l1})+(\mathrm{l2})+(\mathrm{l3})+(\mathrm{l4})+(\mathrm{l5})\notag\\
 \begin{split}
 =&\sum_{k=1}^{m+1}(-1)^{k+m}(-1)^{|\x_k|(|\x_{k+1}|+\cdots+|\x_{m+1}|)}\\
 &\relphantom{+\sum_{k=1}^{m+1}}\cdot(f(\alpha(\x_1),\cdots,\widehat{\alpha(\x_k)},\cdots,\alpha(\x_{m+1}),~~)\cdot[\x_k,\x_{m+2}]_{\alpha})\bullet_{\alpha} \alpha^{m+1}(z)\end{split}\tag{l1}\\
\begin{split}&+\sum_{i=1}^{m+1} (-1)^{i+m}(-1)^{|\x_i|(|\x_{i+1}|+\cdots+|\x_{m+2}|)}\\
 &\relphantom{-\sum_{i=1}^{n-1}}\cdot(f(\alpha(\x_1),\cdots,\widehat{\alpha(\x_i)},\cdots,\alpha(\x_{m+1}),~~)\cdot\alpha(\x_{m+2}))
 \bullet_{\alpha}\alpha^{m}(\x_i\cdot z)\end{split}\tag{l2}\\
 \begin{split}
 &-\sum_{i=1}^{m+1} (-1)^{i+m}(-1)^{|\x_i|(|f|+|\x_{1}|+\cdots+|\x_{i-1}|)}\\
 &\relphantom{-\sum_{i=1}^{n-1}}\cdot\alpha^{m+1}(\x_i)\cdot((f(\x_1,\cdots,\widehat{\x_i},\cdots,\x_{m+1},~~)\cdot \x_{m+2})\bullet_{\alpha} \alpha^{m}(z))
 \end{split}\tag{l3}\\
 \begin{split}
 &-\!\!\sum_{i=1}^{n-1}\sum_{k=1}^{m+1} (-1)^{m+k}(-1)^{(|f|+|\x_{1}|+\cdots+|\x_{{m+1}}|)(|\x_{m+2}^1|+\cdots+|\x_{m+2}^{i-1}|)} (-1)^{|\x_k|(|\x_{k+1}|+\cdots+|\x_{m+1}|)}\\
 &\relphantom{}\cdot[\alpha^{m+1}(\x_{m+2}^1),\cdots,f(\alpha(\x_1),\cdots,\widehat{\alpha(\x_k)},\cdots,\alpha(\x_{m+1}), \x_k\cdot\x_{m+2}^i),\cdots,\\
 &\relphantom{}\quad \quad\alpha^{m+1}(\x_{m+2}^{n-1}),\alpha^{m+1}(z)]
 \end{split}\tag{l4}\\
 \begin{split}
 &+\!\!\sum_{i=1}^{n-1}\sum_{k=1}^{m+1} (-1)^{m+k}(-1)^{(|f|+|\x_{1}|+\!\cdots+|\x_{{m+1}}|)(|\x_{m+2}^1|+\cdots+|\x_{m+2}^{i-1}|)} (-1)^{|\x_k|(|f|+|\x_1|+\cdots+|\x_{k-1}|)}\\
 &\relphantom{}\cdot[\alpha^{m+1}(\x_{m+2}^1),\cdots,\alpha^{m}(\x_k)\cdot f(\x_1,\cdots,\widehat{\x_k},\cdots,\x_{m+1},\x_{m+2}^i),\cdots,\\
 &\relphantom{}\quad \quad\alpha^{m+1}(\x_{m+2}^{n-1}),\alpha^{m+1}(z)]
 \end{split}\tag{l5}\end{align}
 and
 \begin{align*}
&(\mathrm{l1})=\sum_{k=1}^{m+1}(-1)^{k+m}(-1)^{|\x_k|(|\x_{k+1}|+\cdots+|\x_{m+1}|)}\sum_{i=1}^{n-1}(-1)^{|\x_k|(|\x^{1}_{m+2}|+\cdots+|\x^{i-1}_{m+2}|)}\\
 &\quad \quad\cdot(f(\alpha(\x_1),\cdots,\widehat{\alpha(\x_k)},\cdots,\alpha(\x_{m+1}),~~)\cdot(\alpha(\x^{1}_{m+2})\wedge\cdots\wedge\x_k\cdot\x^{i}_{m+2}
 \wedge\cdots\wedge\\
 &\quad \quad\quad\quad\alpha(\x^{n-1}_{m+2}))\bullet_{\alpha} \alpha^{m+1}(z)\\
 &=\sum_{k=1}^{m+1}(-1)^{k+m}(-1)^{|\x_k|(|\x_{k+1}|+\cdots+|\x_{m+1}|)}\sum_{i=1}^{n-1}(-1)^{|\x_k|(|\x^{1}_{m+2}|+\cdots+|\x^{i-1}_{m+2}|)}\\
 &\Big\{\sum_{j<i}(-1)^{(|f|+|\x_{1}|+\cdots+|\widehat{\x_{k}}|+\cdots+|\x_{m+1}|)(|\x^{1}_{m+2}|+\cdots+|\x^{j-1}_{m+2}|)}\\
 &\quad \quad\cdot [\alpha^{m+1}(\x^{1}_{m+2}),\cdots,f(\alpha(\x_1),\cdots,\widehat{\alpha(\x_k)},\cdots,\alpha(\x_{m+1}),\alpha(\x^{j}_{m+2}))
 ,\cdots,\\
 &\quad \quad\quad\quad\alpha^{m}(\x_k\cdot\x^{i}_{m+2}),\cdots,\alpha^{m+1}(\x^{n-1}_{m+2}),\alpha^{m+1}(z)]\tag{l1a}\label{l1a}\\
 &+\sum_{j>i}(-1)^{(|f|+|\x_{1}|+\cdots+|\widehat{\x_{k}}|+\cdots+|\x_{m+1}|)(|\x^{1}_{m+2}|+\cdots+|\x^{j-1}_{m+2}|+|\x_{k}|)}\\
 &\cdot [\alpha^{m+1}(\x^{1}_{m+2}),\cdots,\alpha^{m}(\x_k\cdot\x^{i}_{m+2})
 ,\cdots, f(\alpha(\x_1),\cdots,\widehat{\alpha(\x_k)},\cdots,\alpha(\x_{m+1}),\alpha(\x^{j}_{m+2})),\\
 &\quad \quad\quad\quad \cdots,\alpha^{m+1}(\x^{n-1}_{m+2}),\alpha^{m+1}(z)]\tag{l1b}\label{l1b}\\
 &+(-1)^{(|f|+|\x_{1}|+\cdots+|\widehat{\x_{k}}|+\cdots+|\x_{m+1}|)(|\x^{1}_{m+2}|+\cdots+|\x^{i-1}_{m+2}|)}\\
 &\quad \quad\cdot [\alpha^{m+1}(\x^{1}_{m+2}),\cdots,f(\alpha(\x_1),\cdots,\widehat{\alpha(\x_k)},\cdots,\alpha(\x_{m+1}),\x_k\cdot\x^{i}_{m+2})
 ,\\
 &\quad \quad\quad\quad\cdots,\alpha^{m+1}(\x^{n-1}_{m+2}),\alpha^{m+1}(z)]\Big\}\tag{l1c}\label{l1c}.
 \end{align*}
 Moreover, we have
 \begin{align*}(\mathrm{l3})
&=-\sum_{i=1}^{m+1} (-1)^{i+m}(-1)^{|\x_i|(|f|+|\x_{1}|+\cdots+|\x_{i-1}|)}\\
 &\quad\quad\cdot\sum_{j=1}^{n-1}(-1)^{(|f|+|\x_{1}|+\cdots+|\widehat{\x_{i}}|
+\cdots+|\x_{m+1}|)(|\x^{1}_{m+2}|+\cdots+|\x^{j-1}_{m+2}|)}\\
 &\alpha^{m+1}(\x_i)\cdot[\alpha^{m}(\x^{1}_{m+2}),\cdots, f(\x_1,\cdots,\widehat{\x_i},\cdots,\x_{m+1},\x^{j}_{m+2}),\cdots, \alpha^{m}(\x^{n-1}_{m+2}), \alpha^{m}(z)]\\
 &=-\sum_{i=1}^{m+1} (-1)^{i+m}(-1)^{|\x_i|(|f|+|\x_{1}|+\cdots+|\x_{i-1}|)}\\
 &\quad\cdot\sum_{j=1}^{n-1}(-1)^{(|f|+|\x_{1}|+\cdots+|\widehat{\x_{i}}|
 +\cdots+|\x_{m+1}|)(|\x^{1}_{m+2}|+\cdots+|\x^{j-1}_{m+2}|)}\\
 &\quad\cdot\Big\{\sum_{l<j}(-1)^{|\x_{i}|(|\x^{1}_{m+2}|+\cdots+|\x^{l-1}_{m+2}|)}[\alpha^{m+1}(\x^{1}_{m+2}),
 \cdots,\alpha^{m}(\x_{i})\cdot\alpha^{m}(\x^{l}_{m+2}),\\
 &\quad\quad\cdots,\alpha\circ f(\x_1,\cdots,\widehat{\x_i},\cdots,\x_{m+1},\x^{j}_{m+2})
 ,\cdots,\alpha^{m+1}(\x^{n-1}_{m+2}),\alpha^{m+1}(z)]\tag{l3a}\label{l3a}\\
 &+\sum_{l>j}(-1)^{|\x_{i}|(|f|+|\x_{1}|+\cdots+|\widehat{\x_{i}}|+\cdots+|\x_{m+1}|+|\x^{1}_{m+2}|+\cdots+|\x^{l-1}_{m+2}|)}\\
 &\quad\cdot [\alpha^{m+1}(\x^{1}_{m+2}),\cdots,\alpha\circ f(\x_1,\cdots,\widehat{\x_i},\cdots,\x_{m+1},\x^{j}_{m+2}),\cdots,\alpha^{m}(\x_{i})\cdot\alpha^{m}(\x^{l}_{m+2})\\
 &\quad\quad,\cdots,\alpha^{m+1}(\x^{n-1}_{m+2}),\alpha^{m+1}(z)]\tag{l3b}\label{l3b}\\
 &+(-1)^{|\x_{i}|(|\x^{1}_{m+2}|+\cdots+|\x^{j-1}_{m+2}|)}\cdot [\alpha^{m+1}(\x^{1}_{m+2}),\cdots,\alpha^{m}(\x_{i})\cdot  f(\x_1,\cdots,\widehat{\x_i},\cdots,\\
 &\quad \quad\x_{m+1},\x^{j}_{m+2}),\cdots,\alpha^{m+1}(\x^{n-1}_{m+2}),\alpha^{m+1}(z)]\tag{l3c}\label{l3c}\\
 &+(-1)^{|\x_{i}|(|\x_{m+2}|+|f|+|\x_{1}|+\cdots+|\widehat{\x_{i}}|+\cdots+|\x_{m+1}|)}\cdot [\alpha^{m+1}(\x^{1}_{m+2}),\cdots,\alpha\circ  f(\x_1,\cdots,\widehat{\x_i},\cdots,\\
 &\quad \quad\x_{m+1},\x^{j}_{m+2}),\cdots,\alpha^{m+1}(\x^{n-1}_{m+2}),\alpha^{m}(\x_{i})\cdot\alpha^{m}(z)]\Big\}\tag{l3d}\label{l3d}.
 \end{align*}
 Since $(\mathrm{l4})+(\mathrm{l1c})=0, (\mathrm{l2})+(\mathrm{l3d})=0, (\mathrm{l1a})+(\mathrm{l3b})=0, (\mathrm{l1b})+(\mathrm{l3a})=0,$ one gets $(\mathrm{l1})+(\mathrm{l2})+(\mathrm{l3})+(\mathrm{l4})+(\mathrm{l5})=0.$
 Then $\delta^{m+1}\circ \delta^{m}=0.$ Therefore, we get the following theorem.

\bthm
Let $f\in C^{m}(\g, V)$ be an $m$-cochain. Then $\delta^{m+1}\circ \delta^{m}(f)=0.$
\ethm

\bre The $m$-coboundary operator $\delta^{m}$ as above is a generalization
of the one defined for $n$-ary multiplicative Hom-Nambu-Lie algebras in \cite{FSA} and for first-class $n$-Lie
superalgebras in \cite{YL}. \ere

The map $f\in C^{m}(\g, V)$ is called an $m$-supercocycle if $\delta^{m} f=0$. We denote by $Z^{m}(\g,V)$ the graded subspace spanned by $m$-supercocycles. Since $\delta^{m+1}\circ \delta^{m}(f)=0$ for all $f \in C^{m}(\g, V)$, $\delta^{m-1} C^{m-1}(\g, V)$ is a graded subspace of $Z^{m}(\g,V)$. Therefore we can define a graded cohomology space $H^{m}(\g,V)$ of $\g$ as the graded factor space $Z^{m}(\g,V)/\delta^{m-1} C^{m-1}(\g, V).$
\section{Extensions of $n$-ary multiplicative Hom-Nambu-Lie superalgebras}
\bdefn
Let $(\g_{i}, [\cdot,\cdots,\cdot]_{i},\alpha_{i})(i=1,2,\cdots)$ be a family of $n$-ary multiplicative Hom-Nambu-Lie superalgebras over $\K$. $f_{i}: g_{i}\rightarrow g_{i+1}$ is a morphism of $n$-ary multiplicative Hom-Nambu-Lie superalgebras. The sequence
 $$\xymatrix{\g_{1}\ar[r]^{f_{1}}& \g_{2}\ar[r]^{f_{2}}&\cdots \ar[r]&\g_{i}\ar[r]^{f_{i}}& \g_{i+1}\ar[r]^{f_{i+1}}&\cdots}$$
 is called an exact sequence of $n$-ary multiplicative Hom-Nambu-Lie superalgebras, if it satisfies $\mathrm{Ker}f_{i+1}=f_{i}(g_{i})(i=1,2,\cdots).$
\edefn
\bdefn
Let $(\g, [\cdot,\cdots,\cdot]_{\g},\alpha_{\g}), (\mathfrak{a}, [\cdot,\cdots,\cdot]_{\mathfrak{a}},\alpha_{\mathfrak{a}})$ and
$(\mathfrak{b}, [\cdot,\cdots,\cdot]_{\mathfrak{b}},\alpha_{\mathfrak{b}})$ be $n$-ary multiplicative Hom-Nambu-Lie superalgebras over $\K$. $\g$ is called an extension of $\mathfrak{b}$ by $\mathfrak{a}$ if there is an exact sequence of $n$-ary multiplicative Hom-Nambu-Lie superalgebras:
$$\xymatrix{0\ar[r]& \mathfrak{a}\ar[r]^\iota& \g\ar[r]^\pi& \mathfrak{b}\ar[r]& 0}.$$
\edefn
Let $(\g, [\cdot,\cdots,\cdot]_{\g},\alpha)$ and $(\mathfrak{b}, [\cdot,\cdots,\cdot]_{\mathfrak{b}},\beta)$ be two $n$-ary multiplicative Hom-Nambu-Lie superalgebras over $\K$.
Suppose that $\mathfrak{a}$ is an abelian graded ideal of $\g$, i.e., $\mathfrak{a}$ is a graded ideal such that $[\mathfrak{a},\mathfrak{a},\underbrace{\g,\cdots,\g}_{n-2}]=0$. We consider the case that $\g$ is an extension of $\mathfrak{b}$ by an abelian graded ideal $\mathfrak{a}$ of $\g$. Let $\tau:\mathfrak{b}\rightarrow \g$ be a homogeneous even linear map with $\pi\circ\tau=\id_{\mathfrak{b}}$ and $\alpha\circ \tau=\tau\circ \beta.$
Let $\B=b_1\wedge\cdots\wedge b_{n-1}\in\mathfrak{b}^{\wedge^{n-1}}$ and let $\rho: \mathfrak{b}^{\wedge^{n-1}}\rightarrow \End(\mathfrak{a}), \B\mapsto\tau(\B)=\tau(b_1)\wedge\cdots\wedge\tau(b_{n-1})$. Then $\mathfrak{a}$ becomes a graded $\mathfrak{b}$-module. Let us write $\tau(b)=(0,b)$ and then denote the elements of $\g$ by $(a, b)$ for all $a\in \mathfrak{a}$ and $b\in\mathfrak{b}$. Then, the bracket in $\g$ is defined by
\beq\label{eq:n-def}
[(a_1,b_1),\cdots,(a_n,b_n)]
=\left(\sum_{i=1}^n[\tau(b_1),\cdots,a_i,\cdots,\tau(b_n)]+f(\B,b_n), ~\B\cdot b_n\right),
\eeq
where  $f(\B,b_n)=\tau(\B)\cdot\tau(b_n)-\tau(\B\cdot b_n)$ and  $|(a_i,b_i)|=|a_i|=|b_i|, \forall 1\leq i\leq n$. It is easy to see that $f\in C^1(\mathfrak{b},\mathfrak{a})_{\bar{0}}$.
Let $\A=a_1\wedge\cdots\wedge a_{n-1},(\A,\B)=(a_1,b_1)\wedge\cdots\wedge(a_{n-1},b_{n-1})$ and $(\alpha(\A),\beta(\B))=(\alpha(a_1),\beta(b_1))\wedge\cdots\wedge(\alpha(a_{n-1}),\beta(b_{n-1}))$. Then
\begin{align*}
&(\alpha(\A),\beta(\B))\cdot((\A',\B')\cdot(a_n',b_n'))\\
&-\sum_{i=1}^n(-1)^{|\A|(|a_1'|+\cdots+|a_{i-1}'|)}[(\alpha(a_1'),\beta(b_1')),\cdots,(\A,\B)\cdot(a_i',b_i'),\cdots,(\alpha(a_n'),\beta(b_n'))]\\
=&(\alpha(\A),\beta(\B))\cdot\left(\sum_{i=1}^n[\tau(b_1'),\cdots,a_i',\cdots,\tau(b_n')]+f(\B', b_n'), ~\B'\cdot b_n'\right)\\
&-\sum_{i=1}^n(-1)^{|\A|(|a_1'|+\cdots+|a_{i-1}'|)}\\
&\cdot\left[(\alpha(a_1'),\beta(b_1')),\cdots,
\left(\left(\begin{aligned}
&\sum_{j=1}^{n-1}[\tau(b_1),\cdots, a_j,\cdots,\\
&\quad \tau(b_{n-1}),\tau(b_i')]\\
&+\tau(\B)\cdot a_i'+f(\B,b_i')
\end{aligned}\right)
, ~\B\cdot b_i'\right)
,\cdots,(\alpha(a_n'),\beta(b_n'))\right]\\
=&\left(\left(\begin{aligned}
&\tau(\beta(\B))\cdot\left(\sum_{i=1}^n[\tau(b_1'),\cdots,a_i',\cdots,\tau(b_n')]\right)\\
&+\sum_{j=1}^{n-1}\left[\tau(\beta(b_1)),\cdots, \alpha(a_j),\cdots,\tau(\beta(b_{n-1})),\tau(\B'\cdot b_n')\right]\\
&+\tau(\beta(\B))\cdot f(\B',b_n')+f(\beta(\B), \B'\cdot b_n')
\end{aligned}\right)
, ~\beta(\B)\cdot(\B'\cdot b_n)\right)\\
&-\sum_{i=1}^n(-1)^{|\A|(|a_1'|+\cdots+|a_{i-1}'|)}\\
&\cdot\left(\left(\begin{aligned}
&\sum_{j=1}^{n-1}\Big[\tau(\beta(b_1')),\cdots,[\tau(b_1),\cdots, a_j,\cdots,\\
&\quad\tau(b_{n-1}),\tau(b_i')],\cdots,\tau(\beta(b_n'))\Big]\\
&+[\tau(\beta(b_1')),\cdots,\tau(\B)\cdot a_i',\cdots,\tau(\beta(b_n'))]\\
&+[\tau(\beta(b_1')),\cdots,f(\B,b_i'),\cdots,\tau(\beta(b_n'))]\\
&+\sum_{j\neq i}[\tau(\beta(b_1')),\cdots,\alpha(a^{'}_j),\cdots,\tau(\B\cdot b_i'),\\
& \quad\cdots,\tau(\beta(b_n'))]\\
&+f(\beta(b_1'),\cdots,\B\cdot b_i',\cdots,\beta(b_n'))
\end{aligned}\right)
, ~[\beta(b_1'),\cdots,\B\cdot b_i',\cdots,\beta(b_n')]\right)\\
=&(\delta^{1} f(\B,\B',b_n'), 0).
\end{align*}
Therefore, $f\in Z^1(\mathfrak{b}, \mathfrak{a})_{\bar{0}}$.

Conversely, suppose that an abelian $n$-ary multiplicative Hom-Nambu-Lie superalgebras $\mathfrak{a}$ is a graded $\mathfrak{b}$-module, $\rho(\B)\cdot a:=\tau(\B)\cdot a$, and $f\in Z^1(\mathfrak{b}, \mathfrak{a})_{\bar{0}}$.
Let $\g:=(\mathfrak{a},\mathfrak{b})=\{(x,y)|x\in \mathfrak{a},y\in \mathfrak{b}\},$ $\alpha^{'}:=\alpha+\beta,$ where $(\alpha+\beta)(x,y)=(\alpha(x),\beta(y)), x\in \mathfrak{a}, y\in \mathfrak{b}.$ Then $(\g,\alpha^{'})$ is an $n$-ary multiplicative Hom-Nambu-Lie superalgebra with the bracket defined by (\ref{eq:n-def}). Then we can define an exact sequence
$$\xymatrix{0\ar[r]& \mathfrak{a}\ar[r]^\iota& \g\ar[r]^\pi& \mathfrak{b}\ar[r]& 0},$$
where $\iota(a)=(a,0), \pi(a,b)=b$. Thus $\g$ is an extension of $\mathfrak{b}$ by $\mathfrak{a}$ and $\iota(\mathfrak{a})$ is an abelian graded ideal of $\g$.

Therefore, we get the following theorem.
\bthm
Suppose that $(\mathfrak{a},[\cdot,\cdots,\cdot]_{\mathfrak{a}},\alpha)$ and $(\mathfrak{b},[\cdot,\cdots,\cdot]_{\mathfrak{b}},\beta)$ are two $n$-ary multiplicative Hom-Nambu-Lie superalgebras over $\K$ and $\mathfrak{a}$ is abelian. Then there is a one-to-one correspondence between extensions of $\mathfrak{b}$ by $\mathfrak{a}$ and $Z^1(\mathfrak{b}, \mathfrak{a})_{\bar{0}}$.
\ethm
\section{$T$*-extensions of $n$-ary multiplicative Hom-Nambu-Lie superalgebras}

Let $(\g, [\cdot,\cdots,\cdot],\alpha)$ be an $n$-ary multiplicative Hom-Nambu-Lie superalgebra and $\g^*$ be its dual space.
Since $\g=\g_{\bar{0}}\oplus\g_{\bar{1}}$ and $\g^*=\g^*_{\bar{0}}\oplus\g^*_{\bar{1}}$ are $\Z_2$-graded vector space, the direct sum $\g\oplus\g^*=(\g_{\bar{0}}\oplus\g^*_{\bar{0}})
\oplus(\g_{\bar{1}}\oplus\g^*_{\bar{1}})$ is a $\Z_2$-graded vector space. In the sequel, whenever $x+f\in \g\oplus\g^*$ appears, it means that $x+f$ is homogeneous and $|x+f|=|x|=|f|$.

\blem
Let $\g^*$ be the dual $\Z_2$-graded vector space of an $n$-ary multiplicative Hom-Nambu-Lie superalgebra $(\g, [\cdot,\cdots,\cdot], \alpha)$. Let us consider the even linear map $\ad^*:\g^{\wedge^{n-1}}\rightarrow \End(\g^*)$ defined by
$$\ad^*(\x)(f)(z)=-(-1)^{|\x||f|}f(\ad\x(z)),$$
for all $\x\in \g^{\wedge^{n-1}}, f\in\g^*$ and $z\in\g$. Then $\ad^*$ is a representation of $\g$ on $\g^*$ if and only if the following conditions hold:
\begin{align}\ad(\x)\ad\alpha(\y)-(-1)^{|\x||\y|}\ad(\y)\ad\alpha(\x)=&\alpha \circ \ad[\x,\y]_{\alpha};\tag{i}\\
\begin{split}\ad(x_{1},\cdots,x_{n-2},y_{i})\ad(\alpha(y_{1}),\cdots,
\widehat{\alpha(y_{i})},\cdots,\alpha(y_{n}))
=&(-1)^{(|x_{1}|+\cdots+|x_{n-2}|)(|y_{1}|+\cdots+\widehat{|y_{i}|}+\cdots+|y_{n}|)}\\
\quad\cdot
\Big\{-\ad(\alpha(y_{1}),\cdots,\widehat{\alpha(y_{i})},\cdots,\alpha(y_{n}))&\ad(x_{1},\cdots,x_{n-2},y_{i})\Big\}
\end{split} \tag{ii}
\end{align}for all $i=1,2,\cdots,n.$
We call the representation $\ad^*$ the coadjoint representation of $\g.$
\elem
\bpf

$(\Rightarrow)$ We firstly prove that the necessity holds. Then by the definition of $\ad^*,$ one gets
\begin{align*}\ad^*(\alpha(\x))\ad^*(\y)(f)(z)
=&-(-1)^{|\x|(|\y|+|f|)}\ad^*(\y)(f)(\ad\alpha(\x)(z))\\
=&-(-1)^{|\x|(|\y|+|f|)}(-(-1)^{|\y||f|}f(\ad(\y)\ad\alpha(\x)(z))\\
=&(-1)^{|\x|(|\y|+|f|)+|\y||f|}f(\ad(\y)\ad\alpha(\x)(z))
\end{align*}
and
\begin{align*}(-1)^{|\x||\y|}\ad^*\alpha(\y)\ad^*(\x)(f)(z)
=&(-1)^{|\x||\y|}(-(-1)^{|\y|(|\x|+|f|)}\ad^*(\x)(f)(\ad\alpha(\y)(z))\\
=&-(-1)^{|\y||f|}(-(-1)^{|x||f|}f(\ad(\x)\ad\alpha(\y)(z))\\
=&(-1)^{(|\x|+|\y|)|f|}f(\ad(\x)\ad\alpha(\y)(z)).
\end{align*}
Moreover, we have
\begin{align*}\ad^*([\x,\y]_{\alpha})\circ\nu(f)(z)
=&-(-1)^{(|\x|+|\y|)|f|)}\nu(f)(\ad[\x,\y]_{\alpha}(z))\\
=&-(-1)^{(|\x|+|\y|)|f|}f(\alpha\circ\ad[\x,\y]_{\alpha}(z)).
\end{align*}
By (\ref{eq:rho2}), we have $\ad(\x)\ad\alpha(\y)-(-1)^{|\x||\y|}\ad(\y)\ad\alpha(\x)=\alpha \circ \ad[\x,\y]_{\alpha}.$
\begin{align*}&\ad^*(\alpha(x_{1}),\cdots,\alpha(x_{n-2}),[y_{1},\cdots,y_{n}])\nu(f)(z)\\
=&-(-1)^{(|x_{1}|+\cdots+|x_{n-2}|+|y_{1}|+\cdots+|y_{n}|)|f|} f[\alpha(x_{1}),\cdots,\alpha(x_{n-2}),[y_{1},\cdots,y_{n}],\alpha(z)]\\
=&-(-1)^{(|x_{1}|+\cdots+|x_{n-2}|+|y_{1}|+\cdots+|y_{n}|)|f|}(-(-1)^{|z|(|y_{1}|+\cdots+|y_{n}|)})\\
&\cdot f[\alpha(x_{1}),\cdots,\alpha(x_{n-2}),\alpha(z),[y_{1},\cdots,y_{n}]]\\
=&-(-1)^{(|x_{1}|+\cdots+|x_{n-2}|+|y_{1}|+\cdots+|y_{n}|)|f|}(-(-1)^{|z|(|y_{1}|+\cdots+|y_{n}|)})\\
&\cdot \sum_{i=1}^{n}(-1)^{(|x_{1}|+\cdots+|x_{n-2}|+|z|)(|y_{1}|+\cdots+|y_{i-1}|)}\\
&\cdot f[\alpha(y_{1}),\cdots,\alpha(y_{i-1}),[x_{1},\cdots,x_{n-2},z,y_{i}],\alpha(y_{i+1}),\cdots,\alpha(y_{n})]\\
=&-\sum_{i=1}^{n}(-1)^{n-i}(-1)^{|y_{i}|(|y_{i+1}|+\cdots+|y_{n}|)+(|x_{1}|+\cdots+|x_{n-2}|+|y_{1}|+\cdots+|y_{n}|)|f|}\\
&\cdot(-1)^{(|x_{1}|+\cdots+|x_{n-2}|)(|y_{1}|+\cdots+\widehat{|y_{i}|}+\cdots+|y_{n}|)}  f[\alpha(y_{1}),\cdots,\widehat{\alpha(y_{i})},\cdots,\alpha(y_{n}),[x_{1},\cdots,x_{n-2},y_{i},z]]\\
=&-\sum_{i=1}^{n}(-1)^{n-i}(-1)^{|y_{i}|(|y_{i+1}|+\cdots+|y_{n}|)+(|x_{1}|+\cdots+|x_{n-2}|+|y_{1}|+\cdots+|y_{n}|)|f|}\\
&\cdot(-1)^{(|x_{1}|+\cdots+|x_{n-2}|)(|y_{1}|+\cdots+\widehat{|y_{i}|}+\cdots+|y_{n}|)} \\
&\cdot f(\ad(\alpha(y_{1}),\cdots,\widehat{\alpha(y_{i})},\cdots,\alpha(y_{n}))\ad(x_{1},\cdots,x_{n-2},y_{i})(z))
\end{align*}
and
\begin{align*}&\sum_{i=1}^{n}(-1)^{n-i}(-1)^{(|x_{1}|+\cdots+|x_{n-2}|)(|y_{1}|+\cdots+\widehat{|y_{i}|}+\cdots+|y_{n}|)+|y_{i}|(|y_{i+1}|+\cdots+|y_{n}|)}\\
&\quad\cdot\ad^*(\alpha(y_{1}),\cdots,\widehat{\alpha(y_{i})},\cdots,\alpha(y_{n}))\ad^*(x_{1},\cdots,x_{n-2},y_{i})(f)(z)\\
=&\sum_{i=1}^{n}(-1)^{n-i}(-1)^{(|x_{1}|+\cdots+|x_{n-2}|)(|y_{1}|+\cdots+\widehat{|y_{i}|}+\cdots+|y_{n}|)+|y_{i}|(|y_{i+1}|+\cdots+|y_{n}|)}\\
&\quad\cdot(-(-1)^{(|y_{1}|+\cdots+\widehat{|y_{i}|}+\cdots+|y_{n}|)(|x_{1}|+\cdots+|x_{n-2}|+|y_{i}|+|f|)})\\
&\quad\cdot\ad^*(x_{1},\cdots,x_{n-2},y_{i})(f)(\ad(\alpha(y_{1}),\cdots,\widehat{\alpha(y_{i})},\cdots,\alpha(y_{n}))(z))\\
=&\sum_{i=1}^{n}(-1)^{n-i}(-1)^{|y_{i}|(|y_{1}|+\cdots+|y_{i-1}|)+|f|(|y_{1}|+\cdots+|y_{n}|+|x_{1}|+\cdots+|x_{n-2}|)}\\
&\quad\cdot f(\ad(x_{1},\cdots,x_{n-2},y_{i})\ad(\alpha(y_{1}),\cdots,\widehat{\alpha(y_{i})},\cdots,\alpha(y_{n}))(z)).
\end{align*}
By (\ref{eq:rho3}),  we obtain \begin{align*}
\ad(x_{1},\cdots,x_{n-2},y_{i})\ad(\alpha(y_{1}),\cdots,
\widehat{\alpha(y_{i})},\cdots,\alpha(y_{n}))
=&(-1)^{(|x_{1}|+\cdots+|x_{n-2}|)(|y_{1}|+\cdots+\widehat{|y_{i}|}+\cdots+|y_{n}|)}\\
\quad\cdot
\Big\{-\ad(\alpha(y_{1}),\cdots,\widehat{\alpha(y_{i})},\cdots,\alpha(y_{i}))&\ad(x_{1},\cdots,x_{n-2},y_{i})\Big\}.
\end{align*}

$(\Leftarrow)$ It is easy to see that the sufficiency holds. The proof is complete.
\epf

Let $\theta$ be a homogeneous $n$-linear map from $\g^{\wedge^n}$ into $\g^*$ of degree 0. Now we define a bracket on $\g\oplus\g^*$:
\beq\label{eq:bracketofTextension}\begin{split}
[x_1+f_1,\cdots,x_n+f_n]_{\theta}=&[x_1,\cdots,x_n]_{\g}+\theta(x_1,\cdots,x_n)\\
&+ \sum_{i=1}^n(-1)^{n-i}(-1)^{|x_i|(|x_{i+1}|+\cdots+|x_n|)} \ad^*(x_1,\cdots,\widehat{x_i},\cdots,x_n)\cdot f_i.
\end{split}\eeq

\bthm Let $(\g, [\cdot,\cdots,\cdot],\alpha)$ be an $n$-ary multiplicative Hom-Nambu-Lie superalgebra. Assume that the coadjoint representation exists. Then
$(\g\oplus\g^*, [\cdot,\cdots,\cdot]_{\theta}, \alpha^{'})$ is an $n$-ary multiplicative Hom-Nambu-Lie superalgebra if and only if $\theta\in Z^1(\g, \g^*)_{\bar{0}},$ where $\alpha^{'}(x+f)=\alpha(x)+f\circ\alpha, \forall x\in \g, y\in\g^*.$
\ethm
\bpf
It's clear that $[\cdot,\cdots,\cdot]_{\theta}$ satisfies (\ref{eq:2}) if and only if $\theta\in C^1(\g, \g^*)_{\bar{0}}$. Let $\x+\F=(x_1+f_1)\wedge\cdots\wedge(x_{n-1}+f_{n-1})$ and $\y+\G=(y_1+g_1)\wedge\cdots\wedge(y_{n-1}+g_{n-1})$. Then we have
\begin{align*}
&(\alpha^{'}(\x+\F))\cdot \left((\y+\G)\cdot (y_n+g_n)\right)\\
=&(\alpha(\x)+\F\circ\alpha)\cdot\Big\{
\sum_{i=1}^n(-1)^{n-i}(-1)^{|y_i|(|y_{i+1}|+\cdots+|y_n|)} \ad^*(y_1,\cdots,\widehat{y_i},\cdots,y_n)\cdot g_i\\
&\relphantom{(\x+\F)\cdot\Big\{}+\y\cdot y_n+\theta(\y,y_n)\Big\}\\
=&\alpha(\x)\cdot(\y\cdot y_n)+\theta(\alpha(\x),\y\cdot y_n)+\ad^*(\alpha(\x))\cdot\theta(\y,y_n)\\
&+\sum_{j=1}^{n-1}(-1)^{n-j}(-1)^{|x_j|(|x_{j+1}|+\cdots+|x_{n-1}|+|\y|+|y_n|)} \ad^*(\alpha(x_1),\cdots,\widehat{\alpha(x_j)},\cdots,\alpha(x_{n-1}),\y\cdot y_n)\\
&\cdot(f_j\circ\alpha)+\sum_{i=1}^n(-1)^{n-i}(-1)^{|y_i|(|y_{i+1}|+\cdots+|y_n|)} \ad^*(\alpha(\x))\cdot(\ad^*(y_1,\cdots,\widehat{y_i},\cdots,y_n)\cdot g_i)
\end{align*}
and
\begin{align*}
&\sum_{i=1}^n(-1)^{|\x|(|y_1|+\cdots+|y_{i-1}|)}[\alpha(y_1)+(g_1\circ\alpha),\cdots,(\x+\F)\cdot (y_i+g_i),\cdots,\alpha(y_n)+(g_n\circ\alpha)]_{\theta}\\
=&\sum_{i=1}^n(-1)^{|\x|(|y_1|+\cdots+|y_{i-1}|)}\bigg[\alpha(y_1)+(g_1\circ\alpha),\cdots, \Big\{\x\cdot y_i+\theta(\x,y_i)+\ad^*(\x)\cdot g_i\\
&+\sum_{j=1}^{n-1}(-1)^{n-j}(-1)^{|x_j|(|x_{j+1}|+\cdots+|x_{n-1}|+|y_i|)}
\ad^*(x_1,\!\cdots\!,\widehat{x_j},\!\cdots\!,x_{n-1},y_i)\cdot f_j\Big\},\!\cdots\!,\\
&\quad \quad\alpha(y_n)+(g_n\circ\alpha)\bigg]_\theta\\
=&\sum_{i=1}^n(-1)^{|\x|(|y_1|+\cdots+|y_{i-1}|)}\bigg\{[\alpha(y_1),\cdots,\x\cdot y_i,\cdots,\alpha(y_n)]+\theta(\alpha(y_1),\cdots,\x\cdot y_i,\cdots,\\
&\quad \quad\alpha(y_n))\\
&\relphantom{+}+\sum_{k<i}(-1)^{n-k}(-1)^{|y_k|(|y_{k+1}|+\cdots+|y_n|+|\x|)}
\ad^*(\alpha(y_1),\cdots,\widehat{\alpha(y_k)},\cdots,\x\cdot y_i,\cdots,\\
&\quad \quad\alpha(y_n))\cdot (g_k\circ\alpha)\\
&\relphantom{+}+\sum_{i<k}(-1)^{n-k}(-1)^{|y_k|(|y_{k+1}|+\cdots+|y_n|)}
\ad^*(\alpha(y_1),\cdots,\x\cdot y_i,\cdots,\widehat{\alpha(y_k)},\cdots,\\
&\quad \quad\alpha(y_n))\cdot (g_k\circ\alpha)\\
&\relphantom{+}+(-1)^{n-i}(-1)^{(|\x|+|y_i|)(|y_{i+1}|+\cdots+|y_n|)}
\ad^*(\alpha(y_1),\cdots,\widehat{\alpha(y_i)},\cdots,\alpha(y_n))\cdot \Big\{
\theta(\x,y_i)+\\
&\ad^*(\x)\cdot g_i+\sum_{j=1}^{n-1}(-1)^{n-j}(-1)^{|x_j|(|x_{j+1}|+\cdots+|x_{n-1}|+|y_i|)}
\ad^*(x_1,\cdots,\widehat{x_j},\cdots,x_{n-1},y_i)\cdot f_j\Big\}\bigg\}.
\end{align*}
Since $[ ,\cdots, ]_{\g}$ satisfies (\ref{eq:31}) and $\ad^*(\x)$ satisfies (\ref{eq:rho3}), it can be concluded that $[ ,\cdots, ]_{\theta}$ satisfies (\ref{eq:31}) if and only if
\begin{align*}
0=&\theta(\alpha(\x),\y\cdot y_n)+\ad^*(\alpha(\x))\cdot\theta(\y,y_n)
  -\sum_{i=1}^n(-1)^{|\x|(|y_1|+\cdots+|y_{i-1}|)}\\
  &\quad\cdot\theta(\alpha(y_1),\cdots,\x\cdot y_i,\cdots,\alpha(y_n))-\sum_{i=1}^n(-1)^{|\x|(|y_1|+\cdots+|y_{i-1}|)}\\
  &\quad\cdot(-1)^{n-i}(-1)^{(|\x|+|y_i|)(|y_{i+1}|+\cdots+|y_n|)}\cdot \ad^*(\alpha(y_1),\cdots,\widehat{\alpha(y_i)},\cdots,\alpha(y_n))\cdot\theta(\x,y_i)\\
 =&\delta\theta(\x,\y,y_n),
\end{align*}
i.e., $\theta\in Z^1(\g, \g^*)_{\bar{0}}$.
\epf

\bdefn
Let $(\g,[\cdot,\cdots,\cdot],\alpha)$ be an $n$-ary multiplicative Hom-Nambu-Lie superalgebra. A bilinear form $\langle ,\rangle_{\g}$ on $\g$ is said to be nondegenerate if
$$\g^\perp=\{x\in \g|\langle x,y\rangle_{\g}=0, \forall y\in \g\}=0;$$
invariant if
$$\langle[x_1, \!\cdots\!, x_{n-1},y]_{\g},z\rangle_{\g}=-(-1)^{(|x_1|+\cdots+|x_{n-1}|)|y|}\langle y,[x_1, \!\cdots\!, x_{n-1},z]_{\g}\rangle_{\g}, \forall x_1, \!\cdots\!, x_{n-1}, y, z\in \g;$$
supersymmetric if
$$\langle x,y\rangle_{\g}=(-1)^{|x||y|}\langle y,x\rangle_{\g};$$
consistent if
$$\langle x,y\rangle_{\g}=0, \forall x, y\in \g, |x|\neq|y|;$$
$\alpha$ is called $\langle ,\rangle_{\g}$-symmetric, if
$$\langle \alpha(x),y\rangle_{\g}=\langle \alpha(y),x\rangle_{\g}, \forall x,y\in \g;$$
a subspace $I$ of $L$ is called isotropic if
$I\subseteq I^{\bot}.$

In this section, we only consider consistent bilinear forms. If $\g$ admits a nondegenerate invariant supersymmetric bilinear form $\langle ,\rangle_{\g}$ such that $\alpha$ is $\langle, \rangle_{\g}$-symmetric, then we call $(\g, [\cdot,\cdots,\cdot]_{\g}, \alpha, \langle ,\rangle_{\g})$ a metric $n$-ary multiplicative Hom-Nambu-Lie superalgebra.
In particular, a metric vector space is a pair $(V,\alpha)$ consisting of a $\Z_2$-graded vector space $V=V_{\bar{0}}\oplus V_{\bar{1}}$ and an endmorphism $\alpha$ of $V$ admitting a nondegenerate invariant supersymmetric bilinear form $\langle ,\rangle_{\g}$ such that $\alpha$ is $\langle ,\rangle_{\g}$-symmetric.
\edefn

\blem
Define a bilinear form $\langle , \rangle_{\theta}:(\g\oplus\g^*)\times (\g\oplus\g^*)\rightarrow \K$ by
$$\langle x+f, y+g \rangle_{\theta}=f(y)+(-1)^{|x||y|}g(x).$$ Then $\langle y+g, x+f \rangle_{\theta}=(-1)^{|x||y|}\langle x+f, y+g \rangle_{\theta}, \langle , \rangle_{\theta}$ is nondegenerate and $\alpha^{'}$ is $\langle, \rangle_{\theta}$-symmetric, where $\alpha^{'}(x+f)=\alpha(x)+f\circ\alpha,x\in \g, f\in \g^*.$ Moreover, $(\g\oplus\g^*, [\cdot,\cdots,\cdot]_{\theta},\alpha^{'}, \\\langle , \rangle_{\theta})$
is metric if and only if the following identity holds:
\beq\label{eq:supercyclic} \theta(\x,y)(z)+(-1)^{|y||z|}\theta(\x,z)(y)=0.\eeq
\elem
\bpf If $x+f$ is orthogonal to all elements of $\g\oplus\g^*,$ then for arbitrary element $y+g\in \g\oplus\g^*,$  we have $f(y)=0$ and $(-1)^{|x||y|}g(x)=0,$
which implies that $x=0$ and $f=0,$ so $\langle , \rangle_{\theta}$ is nondegenerate. Moreover, we have
\begin{align*}
\langle y+g, x+f \rangle_{\theta}=&g(x)+(-1)^{|y||x|}f(y)\\
  &=(-1)^{|x||y|}(f(y)+(-1)^{|x||y|}g(x))\\
  &=(-1)^{|x||y|}\langle x+f, y+g \rangle_{\theta}.
\end{align*}
In addition, one gets
\begin{align*}
\langle \alpha^{'}(x+f), y+g \rangle_{\theta}=&\langle \alpha(x)+f\circ\alpha, y+g \rangle_{\theta}\\
  =&f\circ\alpha(y)+(-1)^{|x||y|}g(\alpha(x))
\end{align*}
and
\begin{align*}
\langle x+f, \alpha^{'}(y+g) \rangle_{\theta}=&\langle x+f, \alpha(y)+g\circ\alpha \rangle_{\theta}\\
  =&f\alpha(y)+(-1)^{|x||y|}g\circ\alpha(x).
\end{align*}
Hence, $\langle \alpha^{'}(x+f), y+g \rangle_{\theta}=\langle x+f, \alpha^{'}(y+g) \rangle_{\theta}.$

Furthermore,
$(\g\oplus\g^*, \langle , \rangle_{\theta})$ is metric if and only if
\begin{align*}
0=&\langle(\x+\F)\cdot(y+g),z+h\rangle_{\theta}
+(-1)^{|\x||y|}\langle y+g,(\x+\F)\cdot(z+h)\rangle_{\theta}\\
=&\langle\x\cdot y+\theta(\x,y)+\ad^*(\x)\cdot g, z+h\rangle_{\theta}\\ &+\left\langle\sum_{i=1}^{n-1}(-1)^{n-i}(-1)^{|x_i|(|x_{i+1}|+\cdots+|x_{n-1}|+|y|)}\ad^*\!(x_1,\cdots,\widehat{x_i},\cdots,x_{n-1},y)\cdot f_i, z+h\right\rangle_{\theta}\\
&+(-1)^{|\x||y|}\left\langle y+g, \x\cdot z+\theta(\x, z)+ \ad^*(\x)\cdot h\right\rangle_{\theta}\\
&+(-1)^{|\x||y|}\!\left\langle \!y+g, \!\sum_{i=1}^{n-1}(-1)^{n-i}(-1)^{|x_i|(|x_{i+1}|+\cdots+|x_{n-1}|+|z|)}\ad^*(x_1,\!\cdots\!,\widehat{x_i},\!\cdots\!,x_{n-1},z)\!\cdot\! f_i\!\right\rangle_{\theta}\\
=&\theta(\x,y)(z)+(-1)^{|y||z|}\theta(\x,z)(y),
\end{align*}
i.e., (\ref{eq:supercyclic}) holds.\epf

Now we give the definition of $T^*$-extensions.
\bdefn
For a 1-supercocycle $\theta$ satisfying (\ref{eq:supercyclic}) we shall call the metric $n$-ary multiplicative Hom-Nambu-Lie superalgebra $(\g\oplus \g^{*},[\cdot,\cdots,\cdot]_{\theta},\alpha^{'},\langle , \rangle_{\theta})$ the $T^*$-extension of $(\g,[\cdot,\cdots,\cdot],\alpha)$ (by $\theta$) and denote it by $(T_\theta^*\g,[\cdot,\cdots,\cdot]_{\theta},\alpha^{'}).$
\edefn
\bthm
Let $(\g,[\cdot,\cdots,\cdot],\alpha)$ be an $n$-ary multiplicative Hom-Nambu-Lie superalgebra over a field $\K$. Let
\Beq \g^{(0)}=\g, \g^{(m+1)}=[\g^{(m)},\cdots,\g^{(m)}]_{\g} \text{~~and~~} \g^{1}=\g, \g^{m+1}=[\g^{m},\g,\cdots,\g]_{\g}, \forall m\geq0. \Eeq
$\g$ is called solvable (nilpotent) of length $k$ if and only if there is a smallest integer $k$ such that $\g^{(k)}=0$ ($\g^{k}=0$). Then
\begin{enumerate}[(1)]
   \item  If $\g$ is solvable of length $k$, then $T^{*}_{\theta}\g$ is solvable of length $k$ or $k+1$.
   \item  If $\g$ is nilpotent of length $k$, then $T^{*}_{\theta}\g$ is nilpotent of length at least $k$ and at most $2k-1$. In particular, the nilpotent length of $T^{*}_{0}\g$ is $k$.
   \item  If $\g$ can be decomposed into a direct sum of two Hom-ideals of $\g$, then $T^{*}_{0}\g$ can be too.
\end{enumerate}
\ethm
\bpf
(1) Suppose that $\g$ is solvable of length $k$. Since
$(T^{*}_{\theta}\g)^{(m)}/\g^{*}\cong \g^{(m)}$ and $\g^{(k)}=0$, we have
 $(T^{*}_{\theta}\g)^{(k)}\subseteq \g^{*}$, which implies $(T^{*}_{\theta}\g)^{(k+1)}=0$ because $\g^{*}$ is abelian, and it follows that $T^{*}_{\theta}\g$ is
solvable of length $k$ or $k+1$.

(2) Suppose that $\g$ is nilpotent of length $k$. Since $(T^{*}_{\theta}\g)^{m}/\g^{*}\cong \g^{m}$ and $\g^{k}=0$, we have
$(T^{*}_{\theta}\g)^{k}\subseteq \g^{*}$. Let $f\in(T^{*}_{\theta}\g)^{k}\subseteq \g^{*}, y\in \g$, $\x_{j}+\F_j=(\x_{j}^1+\F_j^1)\wedge\cdots\wedge(\x_{j}^{n-1}+\F_j^{n-1})\in (T^{*}_{\theta}\g)^{\wedge^{n-1}}$, $j=1,\cdots,k-1$. Then
$$\left((\x_{1}+\F_1)\cdots(\x_{k-1}+\F_{k-1})\cdot f\right)(y)
=(\ad^*(\x_1)\cdots \ad^*(\x_{k-1})\cdot f)(y)\in f(\g^k)=0.$$
This proves that $(T^{*}_{\theta}\g)^{2k-1}=0$. Hence $T^{*}_{\theta}\g$ is nilpotent of length at least $k$ and at most $2k-1$.

Now consider the case of trivial $T^*$-extension $T^{*}_{0}\g$ of $\g$. Note that
\begin{align*}
&(\x_{1}+\F_1)\cdots(\x_{k-1}+\F_{k-1})\cdot (y+g)\\
=&\ad(\x_1)\cdots \ad(\x_{k-1})\cdot y+\ad^*(\x_1)\cdots \ad^*(\x_{k-1})\cdot g\\
&+\sum_{j=1}^{k-1}\sum_{i=1}^{n-1}
(-1)^{n-i}(-1)^{|\x_j^i|(|\x_j^{i+1}|+\cdots+|\x_j^{n-1}|+|y|+|\x_{j+1}|+\cdots+|\x_{k-1}|)}\\
&\relphantom{+}\cdot\ad^*(\x_1)\cdots \ad^*(\x_{j-1})\ad^*(\x_j^1,\cdots,\widehat{\x_j^i},\cdots,\x_j^{n-1}, \ad(\x_{j+1})\cdots\ad(\x_{k-1})\cdot y)\cdot\F_j^i\\
=&0.
\end{align*}
Then $(T^{*}_{\theta}\g)^k=0$, as required.

(3) Suppose that $0\neq \g=I\oplus J$,  where $I$ and $J$ are two nonzero Hom-ideals of $\g$. Let $I^{*}=\{f\in\g^*| f(J)=0\}$ and $J^{*}=\{f\in\g^*| f(I)=0\}$. Then $I^{*}$(resp. $J^{*}$) can canonically be identified with the dual space of $I$(resp. $J$) and $\g^*\cong I^*\oplus J^*$.

Note that
\begin{align*}
[T^{*}_{0}I,T^{*}_{0}\g,\cdots,T^{*}_{0}\g]_0=&[I\oplus I^*,\g\oplus \g^*,\cdots,\g\oplus \g^*]_0\\                                             =&[I,\g,\cdots,\g]_{\g}+[I^*,\g,\cdots,\g]_0+[I,\g,\cdots,\g,\g^*]_0\\
\subseteq& I\oplus I^*=T^{*}_{0}I,
\end{align*}
since
\begin{align*}
[I^*,\g,\cdots,\g]_0(J)=&I^*([J,\g,\cdots,\g]_{\g})\subseteq I^*(J)=0
\intertext{and}
[I,\g,\cdots,\g,\g^*]_0(J)=&\g^*([I,J,\g,\cdots,\g]_{\g})=\g^*(0)=0.
\end{align*}
Moreover, for $x+f\in T^{*}_{0}I=I\oplus I^{*},$ we have $\alpha^{'}(x+f)=\alpha(x)+f\circ \alpha\in I\oplus I^{*}$ since $f\circ \alpha\in \g^*$ and $f\circ \alpha(J)\in f(J)=0,$ that is, $\alpha^{'}(T^{*}_{0}I)\subseteq T^{*}_{0}I.$
Then $T^{*}_{0}I$ is a Hom-ideal of $T^{*}_{0}\g$ and so is $T^{*}_{0}J$ in the same way. Hence $T^{*}_{0}\g$ can be decomposed into the direct sum $T^{*}_{0}I\oplus T^{*}_{0}J$ of two nonzero Hom-ideals of $T^{*}_{0}\g$.
\epf

\blem\label{lemma3.1}
Let $(\g,,[\cdot,\cdots,\cdot],\alpha,\langle , \rangle_{\theta})$ be a metric $n$-ary multiplicative Hom-Nambu-Lie superalgebra of even dimension $m$ over a field $\K$ and $I$ be an isotropic $m/2$-dimensional Hom-ideal of $\g$. Then $I$ is abelian.
\elem
\bpf
Since dim$I$+dim$I^{\bot}=m/2+\dim I^{\bot}=m$ and $I\subseteq I^{\bot}$, we have $I=I^{\bot}$.

By $I$ is a Hom-ideal of $\g$, one gets
$$\langle\g, [\g,\cdots,\g,I,I]_{\g}\rangle_{\theta}=\langle[\g,\cdots,\g,I]_{\g},I\rangle_{\theta}\subseteq \langle I,I\rangle_{\theta}=0,$$
which implies $[\g,\cdots,\g,I,I]_{\g}\subseteq \g^{\bot}=0$.
\epf
\bdefn
Let $(\g, [\cdot,\cdots,\cdot]_{\g},\alpha)$ and $(\g^{'}, [\cdot,\cdots,\cdot]_{\g^{'}},\beta)$ be two $n$-ary Hom-Nambu-Lie superalgebras. A linear isomorphism map $\phi: \g\rightarrow \g^{'}$ is called an isomorphism of $n$-ary  Hom-Nambu-Lie superalgebras, if
$$\phi\circ \alpha=\beta\circ \phi;$$
$$\phi[x_{1},\cdots,x_{n}]_{\g}=[\phi(x_{1}),\cdots,\phi(x_{n})]_{\g^{'}}, \forall x_{1},x_{2},\cdots,x_{n}\in \g.$$
\edefn
\bdefn
Two metric $n$-ary multiplicative Hom-Nambu-Lie superalgebras $(\g, [\cdot,\cdots,\\\cdot]_{\g}, \alpha, \langle , \rangle_{\g})$ and  $(\g^{'}, [\cdot,\cdots,\cdot]_{\g^{'}},\beta, \langle , \rangle_{\g^{'}})$ is said to be isometric if there exists an $n$-ary multiplicative Hom-Nambu-Lie superalgebra isomorphism $\phi :\g\rightarrow \g^{'}$ such that $\langle x,y\rangle_{\g}=\langle \phi(x),\phi(y)\rangle_{\g^{'}}, \forall x, y\in \g.$
\edefn

\bthm\label{theorem3.1}
Let $(\g,[\cdot,\cdots,\cdot]_{\g},\beta,\langle , \rangle_{\g})$ be a metric $n$-ary multiplicative Hom-Nambu-Lie superalgebra of dimension $m$ over a field $\K$ of characteristic not 2. Suppose that  $(T_{\theta}^{*}\g_1,[\cdot,\cdots,\cdot]_{\theta},\alpha^{'},\langle , \rangle_{\theta})$ is a $T^*$-extension of $(\g_{1},[\cdot,\cdots,\cdot]_{\g_{1}},\alpha).$ Then $(\g,[\cdot,\cdots,\cdot]_{\g},\beta,\langle , \rangle_{\g})$ is isometric to $(T_{\theta}^{*}\g_1,[\cdot,\cdots,\cdot]_{\theta},\alpha^{'},\langle , \rangle_{\theta})$ if and only if $m$ is even and $\g$ contains an isotropic Hom-ideal $I$ of dimension $m/2$. In particular, $\g_1\cong \g/I$.
\ethm
\bpf
($\Longrightarrow$) Since dim$\g_1$ = dim$\g_1^{*}$, dim$\g$ = dim$T^{*}_{\theta}\g_1=m$ is even. Moreover, $\alpha^{'}(f)=f\circ \alpha\in \g_{1}^*$ for all $f\in \g_{1}^*.$ It is clear that $\g_1^{*}$ is a Hom-ideal of  dimension $m/2$ and by the definition of $\langle , \rangle_{\theta}$, we have $\langle \g_1^*,\g_1^* \rangle_{\theta}=0$, i.e., $\g_1^*$ is isotropic.

($\Longleftarrow$) Suppose that $I$ is an $m/2$-dimensional isotropic graded ideal of $\g$. By Lemma \ref{lemma3.1},  $I$ is abelian. Let $\g_1=\g/I$ and $\pi: \g \rightarrow \g_1$ be the canonical projection. Since $\ch \K\neq2$,  we can choose a complement graded subspace $\g_{0}\subseteq\g$ such that $\g=\g_0\dotplus I$ and $\g_0\subseteq \g_0^{\bot}$. Then $\g_0^{\bot}=\g_0$ since dim$\g_0=m/2$.

Denote by $p_{0}$ (resp. $p_1$) the projection $\g \rightarrow \g_0$ (resp. $\g\rightarrow I$) and let $f^*_1$ denote the homogeneous linear map $I \rightarrow \g_1^{*}: z \mapsto f^*_1(z)$, where $f^*_1(z)(\pi(x)):= \langle z,x\rangle_{\g}, \forall x\in \g, \forall z\in I$.

If $\pi(x)=\pi(y)$, then $x-y\in I$, hence $\langle z,x-y\rangle_{\g}\in \langle z,I\rangle_{\g}=0$ and so $\langle z,x\rangle_{\g}=\langle z,y\rangle_{\g}$, which implies $f^*_1$ is well-defined.  Moreover, $f^*_1$ is bijective and $|f^*_1(z)|=|z|$ for all $z\in I$.

In addition, $f^*_1$ has the following property:
\beq\label{eq:f1*&ad*}\begin{split}
&f^*_1([x_1,\cdots,z_k,\cdots,x_n]_{\g})(\pi(y))\\
=&(-1)^{n-k}(-1)^{|z_k|(|x_{k+1}|+\cdots+|x_n|)}\ad^*(\pi(x_1),\cdots,\widehat{\pi(x_k)},\cdots,\pi(x_n))\cdot f^*_1(z_k)(\pi(y)),
\end{split}\eeq
where $x_1,\cdots,x_{k-1},x_{k+1},\cdots,x_n\in \g$, $z_k\in I$.

Define a homogeneous $n$-linear map
\begin{eqnarray*}
\theta:~~~~~ \g_1\times\cdots\times \g_1~~~~&\longrightarrow&\g_1^{*}\\
(\pi(x_1),\cdots,\pi(x_n))&\longmapsto&f^*_1(p_1([x_1,\cdots,x_n]_{\g})),
\end{eqnarray*}
where $x_1,\cdots,x_n\in \g_0.$ Then $\theta$ is well-defined since $\pi|_{\g_0}:\g_0\rightarrow \g_0/I\cong\g/I=\g_1$ is a linear isomorphism and $\theta\in C^1(\g_1,\g_1^*)_{\bar{0}}$.

Now, define the bracket on $\g_1\oplus \g_1^{*}$ by (\ref{eq:bracketofTextension}), then $(\g_1\oplus \g_1^{*},\alpha^{'})$ is a metric $n$-ary multiplicative Hom-Nambu-Lie superalgebra. Let $\varphi$ be a linear map $\g \rightarrow \g_1\oplus \g_1^{*}$ defined by $\varphi(x+z)=\pi(x)+f^*_1(z), \forall x+z\in \g=\g_0\dotplus I. $
 Since $\pi|_{\g_0}$ and $f^*_1$ are linear isomorphisms, $\varphi$ is also a linear isomorphism. Note that
\begin{align*}
&\varphi([x_1+z_1,\cdots,x_n+z_n]_{\g})=\varphi\left([x_1,\cdots,x_n]_{\g}+\sum_{k=1}^n[x_1,\cdots,z_k,\cdots,x_n]_{\g}\right)\\
=&\varphi\left(p_{0}([x_1,\cdots,x_n]_{\g})+p_1([x_1,\cdots,x_n]_{\g})+\sum_{k=1}^n[x_1,\cdots,z_k,\cdots,x_n]_{\g}\right)\\
=&\pi([x_1,\cdots,x_n]_{\g})+f^*_1\left(p_1([x_1,\cdots,x_n]_{\g})+\sum_{k=1}^n[x_1,\cdots,z_k,\cdots,x_n]_{\g}\right)\\
=&[\pi(x_1),\cdots,\pi(x_n)]_{\g_1}+\theta(\pi(x_1),\cdots,\pi(x_n))\\
&+\sum_{k=1}^n(-1)^{n-k}(-1)^{|z_k|(|x_{k+1}|+\cdots+|x_n|)}\ad^*(\pi(x_1),\cdots,\widehat{\pi(x_k)},\cdots,\pi(x_n))\cdot f^*_1(z_k)\\
=&[\pi(x_1)+f^*_1(z_1),\cdots,\pi(x_n)+f^*_1(z_n)]_{\theta}\\
=&[\varphi(x_1+z_1),\cdots,\varphi(x_n+z_n)]_{\theta},
\end{align*}
where we use the definitions of $\varphi$ and $\theta$ and
(\ref{eq:f1*&ad*}). Moreover, $\varphi\circ\alpha=\alpha^{'}\circ\varphi.$ In fact, for $x+z\in \g=\g_0\dotplus I,$ then
\begin{align*}
&\varphi\circ\alpha(x+z)
=\varphi(\alpha(x)+\alpha(z))\\
=&\pi(\alpha(x))+f^*_1(\alpha(z))
\end{align*}
and \begin{align*}
&\alpha^{'}\circ\varphi(x+z)
=\alpha^{'}(\pi(x)+f^*_1(z))\\
=&\alpha(\pi(x))+f^*_1(z)\circ\alpha.
\end{align*}
Moreover, \begin{align*}
&f^*_1(z)\circ\alpha(\pi(x))
=f^*_1(z)\pi(\alpha(x))\\
=&\langle z,\alpha(x)\rangle_{\g}=\langle \alpha(z),x\rangle_{\g}
=f^*_1(\alpha(z))(\pi(x)).
\end{align*}
Therefore, $f^*_1(z)\circ\alpha=f^*_1(\alpha(z)),$ one gets $\varphi\circ\alpha=\alpha^{'}\circ\varphi.$
Then $\varphi$ is an isomorphism of
$n$-ary multiplicative Hom-Nambu-Lie superalgebras, hence $\g_1\oplus \g_1^{*}$ is an $n$-ary multiplicative Hom-Nambu-Lie superalgebra. Furthermore, we have
\begin{align*}
\langle\varphi(x_{0}+z),\varphi(x_{0}'+z')\rangle_{\theta}
&=\langle\pi(x_{0})+f^*_1(z),\pi(x_{0}')+f^*_1(z')\rangle_{\theta}\\
&=f^*_1(z)(\pi(x_{0}'))+(-1)^{|x_0||x_0'|}f^*_1(z')(\pi(x_{0}))\\
&=\langle z,x_{0}'\rangle_{\g}+(-1)^{|x_0||x_0'|}\langle z',x_{0}\rangle_{\g}=\langle x_{0}+z,x_{0}'+z'\rangle_{\g},
\end{align*}
then $\varphi$ is isometric. The relation
\begin{align*}
&\langle[\varphi(x_1+z_1),\cdots,\varphi(x_n+z_n)]_{\theta},\varphi(x_{n+1}+z_{n+1})\rangle_{\theta}\\
=&\langle\varphi([x_1+z_1,\cdots,x_n+z_n]_{\g}),\varphi(x_{n+1}+z_{n+1})\rangle_{\theta}\\
=&\langle[x_1+z_1,\cdots,x_n+z_n]_{\g},x_{n+1}+z_{n+1}\rangle_{\g}\\
=&-(-1)^{(|x_1|+\cdots+|x_{n-1}|)|x_n|}\langle x_n+z_n, [x_1+z_1,\cdots,x_{n-1}+z_{n-1},x_{n+1}+z_{n+1}]_{\g}\rangle_{\g}\\
=&-(-1)^{(|x_1|+\cdots+|x_{n-1}|)|x_n|}\langle \varphi(x_n+z_n), [\varphi(x_1+z_1),\cdots,\varphi(x_{n-1}+z_{n-1}),\varphi(x_{n+1}+z_{n+1})]_{\theta}\rangle_{\theta}.
\end{align*}
For $x+f,y+g\in \g_1\oplus \g_1^{*},$ then there exist $x^{'}+z_{1},y^{'}+z_{2}\in \g_0\dotplus I$ such that $\varphi(x^{'}+z_{1})=x+f$ and $\varphi(y^{'}+z_{2})=y+g.$ Hence, we have
\begin{align*}
&\langle\alpha^{'}(x+f),y+g\rangle_{\theta}=\langle\alpha^{'}(\varphi(x^{'}+z_{1})),\varphi(y^{'}+z_{2})\rangle_{\theta}\\
=&\langle\varphi(\alpha(x^{'}+z_{1})),\varphi(y^{'}+z_{2})\rangle_{\theta}=\langle\alpha(x^{'}+z_{1}),y^{'}+z_{2}\rangle_{\theta}\\
=&\langle\alpha(y^{'}+z_{2}),x^{'}+z_{1}\rangle_{\theta}=\langle\varphi(\alpha(y^{'}+z_{2})),\varphi(x^{'}+z_{1})\rangle_{\theta}\\
=&\langle\alpha^{'}(\varphi(y^{'}+z_{2})),\varphi(x^{'}+z_{1})\rangle_{\theta}=\langle\alpha^{'}(y+g),x+f\rangle_{\theta}.
\end{align*}
Therefore, $(\g_1\oplus \g_1^{*},[\cdot,\cdots,\cdot]_{\theta},\alpha^{'}, \langle ,\rangle_{\theta})$ is a metric $n$-ary multiplicative Hom-Nambu-Lie superalgebra.
In this way, we get a $T^*$-extension $(T^{*}_{\theta}\g_1,[\cdot,\cdots,\cdot]_{\theta},\alpha^{'},\langle ,\rangle_{\theta})$ of $(\g_{1},[\cdot,\cdots,\cdot]_{\g_{1}},\alpha)$ and consequently, $(\g,[\cdot,\cdots,\cdot]_{\g},\beta,\langle ,\rangle_{\g})$ and
$(T^{*}_{\theta}\g_1,[\cdot,\cdots,\cdot]_{\theta},\alpha^{'},\langle ,\rangle_{\theta})$ are isometric as required.
\epf
Suppose that $(\g,[\cdot,\cdots,\cdot],\alpha)$ is an $n$-ary multiplicative Hom-Nambu-Lie superalgebra and
$\theta_{1}$, $\theta_{2}\in Z^1(\g, \g^{*})_{\bar{0}}$ satisfying
(\ref{eq:supercyclic}). $T^{*}_{\theta_{1}}\g$ and
$T^{*}_{\theta_{2}}\g$ are said to be \textit{equivalent} if there
exists an isomorphism of $n$-ary multiplicative Hom-Nambu-Lie superalgebras $\phi:
T^{*}_{\theta_1}\g\rightarrow  T^{*}_{\theta_2}\g$ such that
$\phi|_{\g^*}=\id_{\g^*}$ and the induced map
$\bar{\phi}:T^{*}_{\theta_1}\g/\g^{*}\rightarrow
T^{*}_{\theta_2}\g/\g^{*}$ is the identity, i.e.,
$\phi(x)-x\in\g^*$. Moreover, if $\phi$ is also an isometry, then
$T^{*}_{\theta_1}\g$ and $T^{*}_{\theta_2}\g$ are said to be
\textit{isometrically equivalent}.

\bprop Suppose that $(\g,[\cdot,\cdots,\cdot],\alpha)$ is an
$n$-ary multiplicative Hom-Nambu-Lie superalgebra over a field $\K$ of characteristic
not 2 and $\theta_{1}$, $\theta_{2}\in Z^1(\g, \g^{*})_{\bar{0}}$
satisfying (\ref{eq:supercyclic}). Then we have
\begin{enumerate}[(1)]
   \item  $(T^{*}_{\theta_{1}}\g,[\cdot,\cdots,\cdot]_{\theta_{1}},\alpha^{'},\langle ,\rangle_{\theta_{1}})$ is equivalent to  $(T^{*}_{\theta_{2}}\g,[\cdot,\cdots,\cdot]_{\theta_{2}},\alpha^{'},\langle ,\rangle_{\theta_{2}})$ if and only if $\theta_1-\theta_2:=\delta\theta'\in\delta C^0(\g, \g^{*})_{\bar{0}}$ and $\theta'(x)\alpha=\theta'(\alpha(x))$ for all $x\in \g.$ Moreover,
       \beq\label{eq:induced bilinear form} \langle x, y\rangle_{\theta'}:=\frac{1}{2}\left(\theta'(x)(y)+(-1)^{|x||y|}\theta'(y)(x)\right)\eeq
       becomes a supersymmetric invariant bilinear form on $\g$ and $\alpha$ is $\langle, \rangle_{\theta'}$-symmetric.
   \item  $T_{\theta_{1}}^{*}\g$ is isometrically equivalent to $T_{\theta_{2}}^{*}\g$ if and only if there is $\theta'\in C^0(\g, \g^{*})_{\bar{0}}$ such that $\theta_1-\theta_2=\delta\theta'$ and the bilinear form induced by $\theta'$ in (\ref{eq:induced bilinear form}) vanishes.
\end{enumerate}
\eprop

\bpf
(1) Let $\phi: T_{\theta_{1}}^{*}\g\rightarrow T_{\theta_{2}}^{*}\g$ be an isomorphism of $n$-ary multiplicative Hom-Nambu-Lie superalgebras satisfying $\phi|_{\g^*}=\id_{\g^*}$ and $\phi(x)-x\in \g^*, \forall x\in \g$.
Set $\theta'(x)=\phi(x)-x$. Then $\theta'\in C^0(\g, \g^{*})_{\bar{0}}$ and
\begin{align}
0=&\phi([x_1+f_1,\cdots,x_n+f_n]_{\theta_1})-[\phi(x_1+f_1),\cdots,\phi(x_n+f_n)]_{\theta_2}\notag\\
=&\phi([x_1,\cdots,x_n]_{\g})+\theta_1(x_1,\cdots,x_n)-[x_1+\theta'(x_1)+f_1,\cdots,x_n+\theta'(x_n)+f_n]_{\theta_2}\notag\\
&+\sum_{i=1}^n(-1)^{n-i}(-1)^{|x_i|(|x_{i+1}|+\cdots+|x_n|)}\ad^*(x_1,\cdots,\widehat{x_i},\cdots,x_n)\cdot f_i\notag\\
\begin{split}
=&\theta'([x_1,\cdots,x_n]_{\g})+\theta_1(x_1,\cdots,x_n)-\theta_2(x_1,\cdots,x_n)\\
&-\sum_{i=1}^n(-1)^{n-i}(-1)^{|x_i|(|x_{i+1}|+\cdots+|x_n|)}\ad^*(x_1,\cdots,\widehat{x_i},\cdots,x_n)\cdot \theta'(x_i)
\end{split}\label{eq: theta}\\
=&\theta_1(x_1,\cdots,x_n)-\theta_2(x_1,\cdots,x_n)-\delta\theta'(x_1,\cdots,x_n).\notag
\end{align}
By $\alpha^{'}\phi=\phi\alpha^{'},$ we may obtain $\theta'(x)\alpha=\theta'(\alpha(x))$ for all $x\in \g.$

For the converse, suppose that $\theta'\in C^0(\g, \g^{*})_{\bar{0}}$ satisfies $\theta_1-\theta_2=\delta\theta'$ and $\theta'(x)\alpha=\theta'(\alpha(x))$ for all $x\in \g.$ Let $\phi: T_{\theta_{1}}^{*}\g\rightarrow T_{\theta_{2}}^{*}\g$ be defined by $\phi(x+f)=x+\theta'(x)+f$. Then $\phi|_{\g^*}=\id_{\g^*}$ and $\phi(x)-x\in \g^*, \forall x\in \g.$ Moreover, $\alpha^{'}\phi=\phi\alpha^{'}.$
In fact,
$$\alpha^{'}\phi(x+f)
=\alpha(x+\theta'(x)+f)
=\alpha(x)+\theta'(x)\alpha+f\alpha
$$
and $$
\phi\alpha^{'}(x+f)
=\phi(\alpha(x)+f\alpha)
=\alpha(x)+\theta'(\alpha(x))+f\alpha.
$$
By $\theta'(x)\alpha=\theta'(\alpha(x)),$ one gets $\alpha^{'}\phi=\phi\alpha^{'}.$ Therefore, $\phi$ is an isomorphism of $n$-ary multiplicative Hom-Nambu-Lie superalgebras, that is, $T^{*}_{\theta_1}\g$ is equivalent to  $T^{*}_{\theta_2}\g$.

It's clear that $\langle ,\rangle_{\theta'}$ defined by (\ref{eq:induced bilinear form}) is supersymmetric. Note that
\begin{align*}
&\langle \x\cdot y, z\rangle_{\theta'}+(-1)^{|\x||y|}\langle y, \x\cdot z\rangle_{\theta'}\\
=&\frac{1}{2}\left(\theta'(\x\cdot y)(z)+(-1)^{(|\x|+|y|)|z|}\theta'(z)(\x\cdot y)\right)\\
 &+\frac{1}{2}(-1)^{|\x||y|}\left(\theta'(y)(\x\cdot z)+(-1)^{(|\x|+|z|)|y|}\theta'(\x\cdot z)(y)\right)\\
=&\frac{1}{2}\bigg\{\theta_2(\x,y)(z)-\theta_1(\x,y)(z)+\ad^*(\x)\theta'(y)(z)\\
 &\relphantom{=\frac{1}{2}}+\sum_{i=1}^{n-1}(-1)^{n-i}(-1)^{|x_i|(|x_{i+1}|+\cdots+|x_{n-1}|+|y|)}\ad^*(x_1,\cdots,\widehat{x_i},\cdots,x_{n-1},y)\cdot\theta'(x_i)(z)\bigg\}\\
 &-\frac{1}{2}(-1)^{|y||z|}\ad^*(\x)\cdot \theta'(z)(y)-\frac{1}{2}\ad^*(\x)\cdot \theta'(y)(z)\\
 &+\frac{1}{2}(-1)^{|y||z|}\bigg\{\theta_2(\x,z)(y)-\theta_1(\x,z)(y)+\ad^*(\x)\theta'(z)(y)\\
 &\relphantom{=\frac{1}{2}}+\sum_{i=1}^{n-1}(-1)^{n-i}(-1)^{|x_i|(|x_{i+1}|+\cdots+|x_{n-1}|+|z|)}\ad^*(x_1,\cdots,\widehat{x_i},\cdots,x_{n-1},z)\cdot\theta'(x_i)(y)\bigg\}\\
=&0,
\end{align*}
where we make use of (\ref{eq: theta})=0 and $\theta_1,\theta_2$ satisfying (\ref{eq:supercyclic}).
Then $\langle ,\rangle_{\theta'}$ is invariant.
In addition, \begin{align*}
&\langle \alpha(x), y\rangle_{\theta'}=\frac{1}{2}\left(\theta'(\alpha(x))(y)+(-1)^{|x||y|}\theta'(y)(\alpha(x))\right)\\
=&\frac{1}{2}\left(\theta'(\alpha(y))(x)+(-1)^{|x||y|}\theta'(x)(\alpha(y))\right)
=\langle \alpha(y), x\rangle_{\theta'}
\end{align*}
since $\theta'(x)\alpha=\theta'(\alpha(x))$ for all $x\in \g.$ That is, $\alpha$ is $\langle x, y\rangle_{\theta'}$-symmetric.

(2) Let the isomorphism $\phi$ be defined as in (1). Then for all $x+f, y+g\in T^{*}_{\theta_{1}}\g$, we have
\begin{align*}
&\langle\phi(x+f),\phi(y+g)\rangle_{\theta_2}=\langle x+\theta'(x)+f,y+\theta'(y)+g\rangle_{\theta_2}\\
=&\theta'(x)(y)+f(y)+(-1)^{|x||y|}\theta'(y)(x)+(-1)^{|x||y|}g(x)\\
=&2\langle x, y\rangle_{\theta'}+\langle x+f,y+g\rangle_{\theta_1}.
\end{align*}
Thus $\phi$ is an isometry if and only if $\langle , \rangle_{\theta'}=0$.
\epf
\blem\label{lemma3.3}
Let $(V,\langle,\rangle_V,\alpha)$ be a metric $\Z_2$-graded vector space of dimension $m$ over an algebraically closed field $\K$ of characteristic not 2 and $\g\subseteq gl(V)$ be a Lie superalgebra consisting of nilpotent homogeneous endomorphisms of $V$ such that for each $f\in \g$, the map $f^{+}:V\rightarrow V$ defined by $\langle f^{+}(v),v'\rangle_V=(-1)^{|f||v|}\langle v,f(v')\rangle_V$ is contained in $\g$, too.
Suppose that $W$ is an isotropic graded subspace of $V$ which is stable under $\g$ and $\alpha,$ i.e., $f(W)\subseteq W$ for all $f\in \g$ and $\alpha(W)\subseteq W,$ then $W$ is contained in a maximally isotropic graded subspace $W_{max}$  of  $V$ which is also stable under $\g$ and $\alpha,$ moreover, $\dim W_{max}=[m/2]$. If $m$ is even, then $W_{max}=W_{max}^\bot$. If $m$ is odd, then $W_{max}\subset W_{max}^\bot, \dim  W_{max}^\bot-\dim  W_{max}=1$, and $f(W_{max}^\bot)\subseteq W_{max}$ for all $f\in \g$.
\elem
\bpf
The proof is by induction on $m$. The base step $m = 0$ is
obviously true. For the inductive step, we consider the following two cases.

Case 1: $W \neq 0$ or there is a nonzero $\g$-stable vector $v\in V$(that is, $\g(v)\subseteq\K v$) such that $\langle v,v\rangle_V=0$.

Case 2: $W=0$ and every nonzero $\g$-stable vector $v\in V$ satisfies $\langle v,v\rangle_V\neq 0.$

In the first case $\K v$ is a nonzero isotropic $\g$-stable graded subspace, and $W^{\perp}$ is also $\g$-stable since $\langle w, f(w^\perp)\rangle_V=(-1)^{|f||w|}\langle f^+(w), w^\perp\rangle_V=0$. Now, consider the bilinear form $\langle , \rangle_{V'}$ on the factor graded space $V'=W^{\perp}/W$ defined by $\langle x^\perp+W, y^\perp+W\rangle_{V'}:=\langle x^\perp, y^\perp\rangle_V$, then $V'$ is metric. Denote by $\pi$ the canonical projection $W^{\perp}\rightarrow V'$ and define $f':V'\rightarrow V'$ by $f'(\pi(w^\perp))=\pi(f (w^\perp))$, then $f'$ is well-defined since $W$ and $W^\perp$ are $\g$-stable. Let $\g':= \{f'|f\in\g\}$. Then $\g'$ is a Lie superalgebra. For each $f\in \g,$ there is a positive integer $k$ such that $f^{k}=0$, which implies that $(f')^{k}=0$. Hence $\g'$ also consists of nilpotent homogeneous endomorphisms of $V'$. Note that $\g'$ satisfies the same conditions of $\g$. In fact, let $x^\perp$ and $y^\perp$ be two arbitrary elements in $W^{\perp}$. Then by the definition of $\langle , \rangle_{V'},$ we have
\begin{align*}
& \langle (f')^{+}(\pi(x^\perp)), \pi(y^\perp)\rangle_{V'}=(-1)^{|f^{'}||x^\perp|}\langle \pi(x^\perp),f'(\pi(y^\perp))\rangle_{V'}\\
=&(-1)^{|f||x^\perp|}\langle \pi(x^\perp),\pi(f(y^\perp))\rangle_{V'}
=(-1)^{|f||x^\perp|}\langle x^\perp,f(y^\perp)\rangle_V\\
=&\langle f^{+}(x^\perp),y^\perp\rangle_V=\langle \pi(f^{+}(x^\perp)),\pi(y^\perp)\rangle_{V'}\\
=&\langle (f^{+})'(\pi(x^\perp)),\pi(y^\perp)\rangle_{V'},
\end{align*}
for arbitrary $f\in \g$, which shows that $(f')^{+}=(f^{+})'\in\g'$ for all $f\in \g$.

Since $\dim V'=\dim W^{\perp}-\dim W=\dim V-2\dim W$, we can use the inductive hypothesis to get a maximally isotropic $\g'$-stable subspace $W'_{max}=W_{max}/W$ in $V'$ and $\alpha(W'_{max})\subseteq W'_{max}.$ Clearly, $\dim W'_{max}$ = $[\frac{\dim V'}{2}]$ = $[\frac{n-2\dim W}{2}]$ = $[n/2]-\dim W$.  For all $x^{\perp}, y^{\perp}\in W_{max}$, the relation $\langle x^{\perp},y^{\perp}\rangle_V=\langle \pi(x^{\perp}),\pi(y^{\perp})\rangle_{V'}=0$ implies that $W_{max}$ is isotropic. Note that $\dim W_{max}=\dim W'_{max}+\dim W=[n/2]$, then $W_{max}$ is maximally isotropic. Moreover, for all $f\in \g$ and $w^{\perp}\in W_{max}$, we have $\pi(f(w^{\perp}))=f'(\pi(w^{\perp}))\in W'_{max}$, which implies $f(w^{\perp})\in W_{max}$. It follows that $W_{max}$ is
$\g$-stable and $\alpha(W_{max})\subseteq W_{max}.$  This proves the first assertion of the lemma in this case.

In the second case, by Engel's Theorem of Lie superalgebras, there is a nonzero $\g$-stable vector $v\in V$ such that $f(v)=0$ for all $f\in \g$. Clearly, $\K v$ is a nondegenerate $\g$-stable graded subspace of $V$, then $V=\K v\dotplus (\K v)^{\perp}$ and $(\K v)^{\perp}$ is also $\g$-stable since $\langle f((kv)^\perp),v\rangle_V=(-1)^{|f||v|}\langle (kv)^\perp,f^+(v)\rangle_V=(-1)^{|f||v|}\langle (kv)^\perp,0\rangle_V=0, \forall f\in \g$. Now, if $(\K v)^{\perp}=0$, then $V=\K v$ and $\g(V)=0$, hence $\g=0$ and so 0 is the maximally isotropic $\g$-stable subspace, then the lemma follows. If $(\K v)^{\perp}\neq0$, then again by Engel's Theorem of Lie superalgebras there is a nonzero $\g$-stable vector $w\in(\K v)^{\perp}\subseteq V$ such that $f(w)=0$ for all $f \in \g$. It follows that $\g$ vanishes on the two-dimensional nondegenerate subspace $\K v\dotplus \K w$ of $V$. Without loss of generality, we can assume that $\langle v,v\rangle_V=1=\langle w,w\rangle_V$.  Set $c=\langle v,w\rangle_V$, then it is easy to check that the nonzero vector $v+(-c +\sqrt{c^{2}-1})w$ is isotropic and $\g$-stable. This contradicts the assumption of Case 2.

Therefore, the existence of a maximally isotropic $\g$-stable graded subspace $W_{max}$ containing $W$ is proved. If $m$ is even, then dim$W_{max}$=dim$W_{max}^\bot=m/2$; if $m$ is odd, then dim$W_{max}^\bot=\frac{m+1}{2}$ and dim$W_{max}=\frac{m-1}{2}$. Since $\g'$ is nilpotent, there exists a nonzero $\pi(w^\bot)\in V'$ such that $\g'(\pi(w^\bot))=0$. Note that dim$V'$=1, which implies $\g'(V')=0$, so $\g(W_{max}^\bot)\subseteq W_{max}$.
\epf
\bthm\label{thm:2}
Let $(\g,[\cdot,\cdots,\cdot],\alpha,\langle ,\rangle_{\g})$ be a nilpotent metric $n$-ary multiplicative Hom-Nambu-Lie superalgebra of dimension $m$ over an algebraically closed field $\K$ of characteristic not 2. If $J$ is an isotropic Hom-ideal of $\g$, then $\g$ contains a maximally Hom-ideal $I$ of dimension $[m/2]$  containing $J$. Moreover, if $m$ is even, then $\g$ is isometric to some $T^{*}$-extension of $\g/I$. If $m$ is odd, then $I^\bot$ is abelian and $\g$ is isometric to a nondegenerate graded ideal of codimension 1 in some $T^{*}$-extension of $\g/I$.
\ethm
\bpf
Consider $\ad(\g^{\wedge^{n-1}})=\{\ad\x|\x\in\g^{\wedge^{n-1}}\}$. Then $\ad(\g^{\wedge^{n-1}})$ is a Lie superalgebra. For any $\x\in \g^{\wedge^{n-1}}$, $\ad\x$ is nilpotent since $\g$ is nilpotent. Then the following identity
$$\langle -\ad\x(y), z\rangle_{\g}=(-1)^{|\x||y|}\langle y, \ad\x(z)\rangle_{\g}$$
implies $(\ad\x)^+:=-\ad\x\in\g$. By $J$ is an isotropic graded ideal of $\g,$ then $J$ is an isotropic $\ad(\g^{\wedge^{n-1}})$-stable graded subspace and $\alpha(J)\subseteq J,$ by Lemma $\ref{lemma3.3},$ so there is a maximally isotropic $\ad(\g^{\wedge^{n-1}})$-stable graded subspace $I$ of $\g$ containing $J$ such that $\alpha(I)\subseteq I$ and $\dim I=[m/2],$ $I$ is also an isotropic graded ideal of $\g.$ Moreover, if $m$ is even, then $\g$ is isometric to some $T^{*}$-extension of $\g/I$ by Theorem \ref{theorem3.1}.

If $m$ is odd, then $\dim I^\bot-\dim I=1$ and $\ad(\g^{\wedge^{n-1}})(I^\bot)\subseteq I$ by Lemma \ref{lemma3.3}. Note that
\begin{align*}
Z(I)=&\{x\in\g|[x,I,\g,\cdots,\g]_{\g}=0\}=\{x\in\g|\langle \g, [x,I,\g,\cdots,\g]_{\g}\rangle_{\g}=0\}\\
    =&\{x\in\g|\langle [I,\g,\cdots,\g]_{\g}, x\rangle_{\g}=0\}=[I,\g,\cdots,\g]_{\g}^\bot=\left(\ad(\g^{\wedge^{n-1}})(I)\right)^{\bot},
\end{align*}
which implies that $I^{\bot}\subset \left(\ad(\g^{\wedge^{n-1}})(I^\bot)\right)^{\bot}=Z(I^{\bot})$, hence $I^{\bot}$ is abelian.

Take any nonzero element $a\notin\g,$ we  define $\alpha^{'}$ by \begin{equation}
\alpha^{'}(x)=
\begin{cases}
a& \text{if $x=a$},\\
\alpha(x)& \text{if $x\in \g$}.
\end{cases}
\end{equation}
Then $\K a$ is a 1-dimensional abelian $n$-ary multiplicative Hom-Nambu-Lie superalgebra. Define a bilinear map $\langle ,\rangle_{a}: \K a\times\K a\rightarrow\K$ by $\langle a,a\rangle_{a}=1$. Then $\langle ,\rangle_{a}$ is a nondegenerate supersymmetric invariant bilinear form on $\K a$. Let $\g'=\g\dotplus\K a$. Define
\begin{align*}
[x_1+k_1 a,\cdots,x_n+k_n a]_{\g'}=&[x_1,\cdots,x_n]_{\g}\\
\intertext{and}
\langle x+k_1 a, y+k_2 a\rangle_{\g'}=&\langle x, y\rangle_{\g}+\langle k_1 a, k_2 a\rangle_{a}.
\end{align*}
Then $(\g',[\cdot,\cdots,\cdot]_{\g^{'}},\alpha^{'}, \langle , \rangle_{\g'})$ is a nilpotent metric $n$-ary multiplicative Hom-Nambu-Lie superalgebra since
\begin{align*}
\langle \alpha^{'}(x+a), y+a \rangle_{\g'}
=&\langle \alpha(x)+\alpha(a), y+a \rangle_{\g'}
=\langle \alpha(x), y\rangle_{\g}+\langle \alpha(a), a\rangle_{a}\\
=&\langle \alpha(y), x\rangle_{\g}+\langle \alpha(a), a\rangle_{a}
=\langle \alpha^{'}(y+a), x+a \rangle_{\g'}
\end{align*}
for all $x,y\in \g$
and $\g$ is a nondegenerate Hom-ideal of codimension 1 of $(\g',\alpha^{'}).$ Since $I^\bot$ is not isotropic and $\K$ is algebraically closed there exists $z\in I^{\bot}$ and $\langle z, z\rangle_{\g}=-1.$ In addition, we have $\alpha(I^\bot)\subseteq I^\bot$ since $\langle \alpha^{'}(v^\bot), v \rangle_{\g}=\langle \alpha(v), v^\bot \rangle_{\g}=0$ for $v\in I$ and $v^\bot\in I^\bot.$ Let $b=a+z$ and $I'=I\dotplus\K b$. Then $I'$ is an $(m+1)/2$-dimensional isotropic graded ideal of $\g'$.

In fact, for all $x+k_1 a+k_1z, y+k_2 a+k_2z\in I'$,
\begin{align*}
\langle x+k_1 a+k_1z, y+k_2 a+k_2z\rangle_{\g'}
=&\langle x+k_1z, y+k_2z\rangle_{\g}+\langle k_1 a, k_2 a\rangle_{a}\\
=&\langle x, y\rangle_{\g}+\langle x, k_2z\rangle_{\g}+\langle k_1z, y\rangle_{\g}+\langle k_1z, k_2z\rangle_{\g}+k_1k_2\\
=&k_1k_2-k_1k_2=0.
\end{align*}
In light of Theorem \ref{theorem3.1}, we conclude that $\g'$ is isometric to some $T^{*}$-extension of $\g'/I'$.

Define $\Phi:\g'\rightarrow \g/I, x+\lambda a\mapsto x-\lambda z+I$.  Then
\begin{align*}
[\Phi(x_1+\lambda_1 a),\cdots,\Phi(x_n+\lambda_n a)]_{\g/I}
=&[x_1-\lambda_1 z+I,\cdots,x_n-\lambda_n z+I]_{\g/I}\\
=&[x_1,\cdots,x_n]_{\g}+I=\Phi([x_1,\cdots,x_n]_{\g})\\
=&\Phi([x_1+\lambda_1 a,\cdots,x_n+\lambda_n a]_{\g'}),
\end{align*}
where we use the fact that $I^\bot$ is abelian and $\ad(\g^{\wedge^{n-1}})(I^\bot)\subseteq I$. Moreover, $\Phi\alpha=\alpha\Phi.$ In fact, for $x+\lambda a\in \g^{'},$ we have \begin{align*}
&\Phi\alpha(x+\lambda a)=\Phi(\alpha(x)+\lambda \alpha(a))
=\Phi(\alpha(x)+\lambda a)
=\alpha(x)-\lambda z+I\\
=&\alpha(x)-\lambda \alpha(z)+I
=\alpha(x-\lambda z+I)=\alpha\Phi(x+\lambda a).
\end{align*} It's clear that $\Phi$ is surjective and Ker$\Phi=I'$, so $\g'/{I'}\cong \g/I$, hence the theorem follows.
\epf

Now we show that there exists an isotropic Hom-ideal in every finite-dimensional metric $n$-ary multiplicative Hom-Nambu-Lie superalgebra and investigate the nilpotent length of $\g/I$.
\bprop\label{prop:2}
Suppose that $(\g,[\cdot,\cdots,\cdot],\alpha,\langle ,\rangle_{\g})$ is a finite-dimensional metric $n$-ary multiplicative Hom-Nambu-Lie superalgebra.
\begin{enumerate}[(1)]
   \item  For any graded subspace $V\subseteq \g$, $C(V):=\{x\in\g|[x,\g,\cdots,\g]_{\g}\subseteq V\}=[\g,\cdots,\g,V^{\bot}]_{\g}^{\bot}$.
   \item  $\g^m=C_m(\g)^{\bot}$, where $C_0(\g)=0, C_{i+1}(\g)=C(C_{i}(\g))$.
   \item  If $\g$ is nilpotent of length $k$, then $\g^i\subseteq C_{k-i}(\g)$.
\end{enumerate}
\eprop
\bpf
The relation
$$\langle C(V),[\g,\cdots,\g,V^{\bot}]_{\g}\rangle_{\g}=\langle [\g,\cdots,\g,C(V)]_{\g}, V^{\bot}\rangle_{\g}\subseteq \langle V, V^{\bot}\rangle_{\g}=0$$
shows that $C(V)\subseteq[\g,\cdots,\g,V^{\bot}]_{\g}^{\bot}$. Notice that
$$\langle [\g,\cdots,\g,[\g,\cdots,\g,V^{\bot}]_{\g}^{\bot}]_{\g}, V^{\bot}\rangle_{\g}=\langle [\g,\cdots,\g,V^{\bot}]_{\g}^{\bot}, [\g,\cdots,\g,V^{\bot}]_{\g}\rangle_{\g}=0,$$
which implies $[\g,\cdots,\g,[\g,\cdots,\g,V^{\bot}]_{\g}^{\bot}]_{\g}\subseteq (V^{\bot})^{\bot}=V$, i.e., $[\g,\cdots,\g,V^{\bot}]_{\g}^{\bot}\subseteq C(V)$. Hence (1) follows.

By induction, (2) and (3) can be proved easily.
\epf

\bthm
Every finite-dimensional nilpotent metric $n$-ary multiplicative Hom-Nambu-Lie superalgebra $(\g,[\cdot,\cdots,\cdot],\alpha,\langle ,\rangle_{\g})$ over  an algebraically closed field of characteristic not 2 such that $\alpha(\g)=\g$ is isometric to (a nondegenerate ideal of codimension 1 of) a $T^*$-extension of a nilpotent $n$-ary multiplicative Hom-Nambu-Lie superalgebra whose nilpotent length is at most a half of the nilpotent length of $\g$.
\ethm
\bpf
Define $J=\sum\limits_{i=0}^{\infty}\g^{i}\cap C_{i}(\g)$. Since $\g$ is nilpotent, the sum is finite. Proposition \ref{prop:2} (2) says $(\g^{i})^{\bot}=C_{i}(\g)$, then $\g^{i}\cap C_{i}(\g)$ is isotropic for all $i\geq 0$.  Since
$$\g^{i}\supseteq \g^{j}\supseteq \g^{j}\cap C_j(g), ~\text{if} ~i<j,$$
we have
$$(\g^{j}\cap C_{j}(\g))^{\bot}\supseteq (\g^{i})^{\bot}=C_{i}(\g)\supseteq C_{i}(\g)\cap \g^{i}, ~\text{if} ~i<j.$$
It follows that
$$\langle\g^{i}\cap C_{i}(\g), \g^{j}\cap C_{j}(\g)\rangle_{\g}=0, ~~\forall i,j\geq0.$$
Therefore $J$ is an isotropic graded ideal of $\g$. Let $k$ denote the nilpotent length of $\g$. Using Proposition \ref{prop:2} (3) we can conclude that $\g^{[(k+1)/2]}\subseteq C_{[(k+1)/2]}(\g)$. This implies that $\g^{[(k+1)/2]}$ is contained in $J$.   By Theorem \ref{thm:2},  there is a maximally isotropic graded ideal $I$ of $\g$ containing $J\supseteq \g^{[(k+1)/2]}$. It means that $\g/I$ has nilpotent length at most $[(k+1)/2]$, and the theorem follows.
\epf

\bre
 Most results concerning $T^*$-extensions in \cite{B,lyc,LZ2,YL} are contained in this section as special cases.
\ere

\end{document}